\pdfoutput=1

\documentclass{thielstyle-ams}    

\usepackage[backend=bibtex, hyperref=true, language=auto, texencoding=auto, bibencoding=auto, isbn=false, doi=true, maxnames=5, minnames=1, sorting=anyt]{biblatex}

\newcommand\code[1]{{\normalfont\fontfamily{cmvtt}\selectfont #1}}

\DeclareMathOperator{\Bl}{Bl}
\DeclareMathOperator{\Frac}{Frac}

\newcommand{\mscr}[1]{\mathscr{#1}}

\usepackage{enumitem}

\newlist{enum_thm}{enumerate}{1}
\setlist[enum_thm,1]
{
    label=(\alph*),
    noitemsep,
    nolistsep,
    align=left,
    labelindent=\parindent,
    leftmargin=*
}

\newlist{enum_bullet}{enumerate}{1}
\setlist[enum_bullet,1]
{
    label=$\bullet$,
    noitemsep,
    nolistsep,
    align=left,
    labelindent=\parindent,
    leftmargin=*
}

\newlist{enum_proof}{enumerate}{1}
\setlist[enum_proof,1]
{
    label=(\alph*),
    noitemsep,
    nolistsep,
    align=left,
    labelindent=\parindent,
    itemindent=*,
    leftmargin=0pt
}

\title{Blocks in flat families of finite-dimensional algebras} 
\author{Ulrich Thiel}
\keywords{}
\bibliography{references}

\makeatletter
\tikzset{
  column sep/.code=\def\pgfmatrixcolumnsep{\pgf@matrix@xscale*(#1)},
  row sep/.code   =\def\pgfmatrixrowsep{\pgf@matrix@yscale*(#1)},
  matrix xscale/.code=%
    \pgfmathsetmacro\pgf@matrix@xscale{\pgf@matrix@xscale*(#1)},
  matrix yscale/.code=%
    \pgfmathsetmacro\pgf@matrix@yscale{\pgf@matrix@yscale*(#1)},
  matrix scale/.style={/tikz/matrix xscale={#1},/tikz/matrix yscale={#1}}}
\def\pgf@matrix@xscale{1}
\def\pgf@matrix@yscale{1}
\makeatother

\begin{document}  

\maketitle
\thispagestyle{empty}

\blfootnote{Date: Dec 6, 2017} \blfootnote{\textsc{ULRICH THIEL}, School of Mathematics and Statistics, University of Sydney, NSW 2006 Australia.  Email: \code{ulrich.thiel@sydney.edu.au}}

\begin{abstract}
We study the behavior of blocks in flat families of finite-di\-men\-sional algebras. In a general setting we construct a finite directed graph encoding a stratification of the base scheme according to the block structures of the fibers. This graph can be explicitly obtained when the central characters of simple modules of the generic fiber are known. We show that the block structure of an arbitrary fiber is completely determined by ``atomic'' block structures living on the components of a Weil divisor. As a byproduct, we deduce that the number of blocks of fibers defines a lower semicontinuous function on the base scheme. We furthermore discuss how to obtain information about the simple modules in the blocks by generalizing and establishing several properties of decomposition matrices by Geck and Rouquier. 
\end{abstract}

\section{Introduction}

It is a classical fact in ring theory that a non-zero noetherian ring $A$ can be decomposed as a direct product $A= \prod_{i=1}^n B_i$ of indecomposable rings $B_i$. Such a decomposition is unique up to permutation and isomorphism of the factors. Let us denote by $\msf{Bl}(A)$ the set of the $B_i$, called the \word{blocks} of $A$. The decomposition of $A$ into blocks induces a decomposition $A\tn{-}\msf{Mod} = \bigoplus_{i=1}^n B_i\tn{-}\msf{Mod}$ of the category of (left) $A$-modules. In particular, a simple $A$-module is a simple $B_i$-module for a unique block $B_i$ and so we get an induced decomposition $\msf{Irr} A = \coprod_{i=1}^n \msf{Irr} B_i$ of the set of simple modules. Let us denote by $\msf{Fam}(A)$ the set of the $\msf{Irr} B_i$, called the \word{families} of $A$. The blocks and families of a ring are important invariants which help to organize and simplify its representation theory. The aim of this paper is to investigate how these invariants vary in a flat family of finite-dimensional algebras. 

More precisely, we consider a \word{finite flat} algebra $A$ over an integral domain $R$, i.e., $A$ is finitely generated and flat as an $R$-module. This yields a family of algebras parametrized by $\Spec(R)$ consisting of the \word{specializations} (or \word{fibers})
\begin{equation} \label{specialization_definition}
A(\fp) \dopgleich \msf{k}(\fp) \otimes_R A \simeq A_\fp/\fp_\fp A_\fp \;,
\end{equation}
where $\msf{k}(\fp) = \msf{Frac}(R/\fp)$ is the residue field of $\fp \in \Spec(R)$ in $R$ and $A_\fp$ is the localization of $A$ in $\fp$. Note that the fiber $A(\fp)$ is a finite-dimensional $\msf{k}(\fp)$-algebra. Now, the primary goal would be to describe for any $\fp$ the blocks of $A(\fp)$, e.g., the number of blocks, and to describe the simple modules in each block, e.g., the number of such modules and their dimensions. 

It is clear that there will be no general theory giving the precise solutions to these problems for arbitrary $A$. For example, we can take the group ring $A = \bbZ \mrm{S}_n$ of the symmetric group. The fibers of $A$ are precisely the group rings $\bbQ \mrm{S}_n$ and $\bbF_p \mrm{S}_n$ for all primes $p$, and the questions above are still unanswered. Nonetheless, and this is the point of this paper, there are some general phenomena, some patterns in the  behavior of blocks and simple modules along the fibers, which are true quite generally.

\subsection{The setting}
We assume that $R$ is noetherian and normal, and that the generic fiber $A^K$ is a \textit{split} $K$-algebra, where $K$ is the fraction field of $R$, i.e., all simple modules of $A^K$ remain simple under field extension. This setting includes many interesting examples in representation theory like Brauer algebras, Hecke algebras, (restricted) rational Cherednik algebras, etc. We note that some results we mention below actually hold more generally and refer to the main body of the paper. 

At the very end we also establish a semicontinuity property of blocks in the (important) case of a non-split generic fiber, see Theorem \ref{non-split-intro}. This then applies also to quantized enveloping algebras of semisimple Lie algebras at roots of unity, enveloping algebras of semisimple Lie algebras in positive characteristic, quantized function algebras of semisimple groups at roots of unity, etc. More generally, this applies to Hopf PI triples as introduced by Brown–Goodearl \cite{Brown-Goodearl-Quantum-Groups} (see also Brown–Gordon \cite{BG-Ramification} and Gordon \cite{Gordon:Representations-of-semisimple-Lie-in-pos}).

\subsection{Block stratification}
Under the assumptions described above, we prove the following theorem (see Corollary \ref{normal_split_corollary}) which is the backbone of this paper:

\begin{thm}
For any $\fp \in \Spec(R)$ the natural map $A_\fp \twoheadrightarrow A(\fp)$ is \textnormal{block bijective}, i.e., it induces a bijection between the block idempotents.
\end{thm}

This allows us to reduce the study of blocks of specializations to blocks of localizations, and this is much simpler from the general perspective. Since a block idempotent of a localization $A_\fp$ splits into a sum of block idempotents of the generic fiber $A^K$, we can view the blocks of $A_\fp$ as being a \textit{partition} of the set of blocks of $A^K$, see \S\ref{blocks_of_localizations} for details. This gives us a direct and natural way of comparing the block structures among the fibers---something which is in general, without the above theorem, not possible. Let
\[
\Bl_A : \Spec(R) \to \msf{Part}( \Bl(A^K) )
\]
be the map just described. We equip the image $\Bl(A)$ of this map with the partial order $\leq$ on partitions, where $\mscr{P}' \leq \mscr{P}$ if the members of $\mscr{P}'$ are unions of members of $\mscr{P}$. We let $\Bl_A^{-1}(\mscr{P})$, respectively $\Bl_A^{-1}(\leq \! \mscr{P})$, be the locus of all $\fp \in \Spec(R)$ such that the block structure of $A_\fp$, and thus of $A(\fp)$ under the above bijection, is equal to a given partition $\mscr{P}$, respectively coarser than $\mscr{P}$. We then obtain as Theorem \ref{block_stratification_thm}:

\begin{thm}
The sets $\Bl_A^{-1}(\leq \! \mscr{P})$ are closed in $\Spec(R)$, the sets $\Bl_A^{-1}(\mscr{P})$ are locally closed in $\Spec(R)$, and $\Spec(R) = \coprod_{\mscr{P}} \Bl_A^{-1}(\mscr{P})$ is a stratification of $\Spec(R)$.
\end{thm}

Denoting by $\bullet$ the generic point of $\Spec(R)$, so that $A_\bullet = A(\bullet) = A^K$, this implies in particular that the set
\[
\msf{BlGen}(A) \dopgleich \Bl_A^{-1}( \Bl_A(\bullet))  = \{ \fp \in \Spec(R) \mid \Bl_A(\fp) = \Bl_A(\bullet) \}
\]
of primes $\fp$ where the block structure of the fiber $A(\fp)$ is equal to the one of the generic fiber $A^K$ is an open (dense) subset of $\Spec(R)$. Hence, the set
\[
\msf{BlEx}(A) \dopgleich \{ \fp \in \Spec(R) \mid \Bl_A(\fp) < \Bl_A(\bullet) \}
\]
of primes where the block structure of the fiber is coarser than the one of the generic fiber is closed. This set has a nice property, see Corollary \ref{blex_is_weil_div}:

\begin{thm}
If $R$ is a Krull domain (e.g., if $R$ is normal), then $\msf{BlEx}(A)$ is a reduced Weil divisor, i.e., it is either empty or pure of codimension one in $\Spec(R)$.
\end{thm}

We thus call $\msf{BlEx}(A)$ the \word{block divisor} of $A$. This is an interesting new \textit{discriminant} of $A$. Let $\msf{At}(A)$ be the set of irreducible components of $\msf{BlEx}(A)$. On any $Z \in \msf{At}(A)$ there is a unique maximal block structure $\msf{Bl}_A(Z)$, namely the one in the generic point. In \S\ref{block_strat} we show that these block structures have an \textit{atomic character}:

\begin{thm}
For $\fp \in \Spec(R)$ we have 
\[
\Bl_A(\fp) = \bigwedge_{ \substack{Z \in \msf{At}(A) \\ \fp \in Z}} \Bl_A(Z)  \;,
\]
where $\wedge$ is the meet of partitions, i.e., the members are the unios of all members with non-empty intersection.
\end{thm}

Hence, once we know $\msf{At}(A)$ and the \word{atomic block structures} $\msf{Bl}_A(Z)$ for $Z \in \msf{At}(A)$, we know the block structure for any $\fp \in \Spec(R)$. By considering sets of the form
\[
\bigcap_{Z \in \mscr{Z}} Z \setminus \bigcup_{Z \notin \mscr{Z}} Z
\]
for subsets $\mscr{Z} \subs \msf{At}(A)$, we obtain a stratification of $\Spec(R)$ refining the one introduced above. We call this the \word{block stratification} of $A$. 

\subsection{Blocks via central characters}
In \S\ref{blocks_via_central_characters} we discuss an approach to explicitly compute the block stratification and the block structures on the strata. This is based on the knowledge of central characters of simple $A^K$-modules. Since $A^K$ splits, each simple $A^K$-module $S$ has a central character $\Omega_S' : \msf{Z}(A) \to R$, the image lying in $R$ since $R$ is normal. In Theorem \ref{maintheorem_blex_descr} we show:

\begin{thm} \label{maintheorem_blex_descr_intro}
Two simple $A^K$-modules $S$ and $T$ lie in the same $A_\fp$-block if and only if $\Omega_S' \equiv \Omega_T' \ \msf{mod} \ \fp$. 
\end{thm}

The key ingredient in the proof is a (rather non-trivial) result by B. Müller stating that the cliques of a noetherian ring, which is a finite module over its center, are fibered over the center. We address this in detail in \S\ref{sec_mueller}.

If $z_1,\ldots,z_n$ is an $R$-algebra generating system of $\msf{Z}(A)$, then $\Omega_S' \equiv \Omega_T' \ \msf{mod} \ \fp$ if and only if $\Omega_S'(z_i) \equiv \Omega_T'(z_i) \ \msf{mod} \ \fp$ for all $i$. Hence, Theorem \ref{maintheorem_blex_descr_intro} gives a computational tool to explicitly determine $\Bl_A(\fp)$ once the central characters of the generic fiber are known. Moreover, it follows from Theorem \ref{maintheorem_blex_descr_intro} that $\msf{At}(A)$ is the set of maximal irreducible components of the zero loci of the sets
\begin{equation}
 \lbrace \Omega_S'(z_i) - \Omega_T'(z_i) \mid i =1,\ldots,n \rbrace
\end{equation}
for $\Omega_S' \neq \Omega_T'$. The atomic block structures can then be determined by the vanishing of the differences $\Omega_S'- \Omega_T'$ on the $Z \in \msf{At}(A)$, and from these we obtain all block structures as described above.

\subsection{An example}

Let us illustrate this with an explicit example. Let $A$ be the generic Brauer algebra for $n=3$ over the polynomial ring $R \dopgleich \bbZ \lbrack \boldsymbol{\delta} \rbrack$, see \cite{Graham-Lehrer-Cellular}. There are four simple $A^K$-modules, labelled by the partitions $(0,(1,1,1))$, $(0,(3))$, $(0,(2,1))$, and $(1,(1))$. We will simply label these by $1,\ldots,4$ from now on. Since $A^K$ is semisimple, we can identify the blocks of $A^K$ with the simple modules of $A^K$, i.e., we can label the blocks by $1,\ldots,4$. We can thus view blocks of specializations of $A$ as partitions of $\{ 1,\ldots,4 \}$ as described above. It is not too difficult to explicitly compute the central characters of the simple $A^K$-modules. From these we deduce that the block structure of the fibers of $A$ over $\bbZ \lbrack \boldsymbol{\delta} \rbrack$ are distributed as in the following graph:
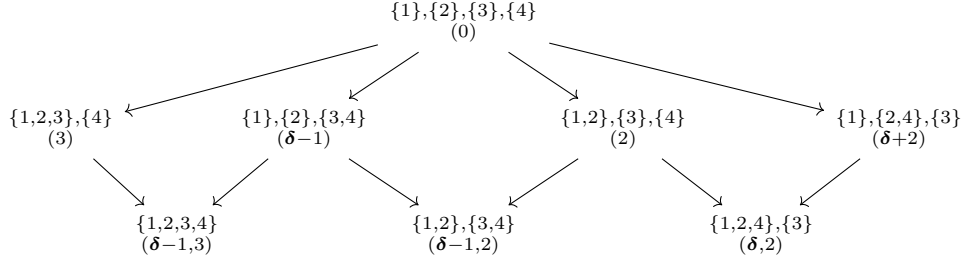
\begin{figure}[htbp]
\begin{tikzcd}[column sep=0pt]
& & & \substack{\lbrace 1 \rbrace, \lbrace 2 \rbrace, \lbrace 3 \rbrace, \lbrace 4 \rbrace \\ (0)} \arrow{dlll} \arrow{dl} \arrow{dr} \arrow{drrr} \\
\substack{\lbrace 1,2,3 \rbrace, \lbrace 4 \rbrace \\ (3)} \arrow{dr} & & \substack{\lbrace 1 \rbrace, \lbrace 2 \rbrace, \lbrace 3,4 \rbrace \\ (\boldsymbol{\delta}-1)} \arrow{dl} \arrow{dr} & & \substack{\lbrace 1,2 \rbrace, \lbrace 3 \rbrace, \lbrace 4 \rbrace \\ (2)} \arrow{dr} \arrow{dl} & & \substack{\lbrace 1 \rbrace, \lbrace 2,4 \rbrace, \lbrace 3 \rbrace \\ (\boldsymbol{\delta}+2)}\arrow{dl} \\
& \substack{\lbrace 1,2,3,4 \rbrace \\ (\boldsymbol{\delta}-1,3)} & & \substack{\lbrace 1,2 \rbrace, \lbrace 3,4 \rbrace \\ (\boldsymbol{\delta}-1,2)} & & \substack{ \lbrace 1,2,4 \rbrace, \lbrace 3 \rbrace \\ (\boldsymbol{\delta},2)} 
\end{tikzcd}
\end{figure}

This graph encodes the block stratification of the two-dimensional base scheme $\Spec(\bbZ \lbrack \boldsymbol{\delta} \rbrack)$, along with the block structures on the strata. We see that $\msf{BlEx}(A)$ has four components of codimension one, the generic points of these components are $3$, $\delta-1$, $2$, and $\delta+2$, respectively. The block structure on any other point $\fp$ is uniquely determined as the meet of the block structures on the components of $\msf{BlEx}(A)$ containing $\fp$. 

We want to point out that it is central for us to work with (affine) schemes. For example, we have one skeleton with generic point $(2)$, i.e., we consider the Brauer algebra in characteristic two. Now, we do not only have the case $\boldsymbol{\delta} \in \lbrace 0,1 \rbrace = \bbF_2$, which is described by the two strata below $(2)$, but we also have a generic characteristic two case, described by the generic point of $\bbF_2 \lbrack \boldsymbol{\delta} \rbrack$, and this is really different from the case of specialized~$\boldsymbol{\delta}$ as we can see from the block structures.

Note that the components of $\msf{BlEx}(A)$ are precisely the parameters where the Brauer algebra is not semisimple anymore (the precise parameters have been determined by Rui \cite{Rui} for all $n \in \bbN$). We show in Lemma \ref{cellular_alg_blgen_decgen} that this is always the case for cellular algebras. 

\subsection{Blocks and decomposition matrices}
 In \S\ref{blocks_and_dec_maps} we address questions about the simple modules in a block. The main tool here are the decomposition matrices introduced by Geck and Rouquier. In Theorem \ref{brauer_rec} we show that they satisfy Brauer reciprocity in a rather general setting in which it was not known to hold before. In \S\ref{brauer_graph} we generalize the concept of Brauer graphs and show how these relate to blocks. 
 
 \subsection{An open problem}
 In \S\ref{preservation} we contrast the preservation of simple modules with the preservation of blocks under specialization, and this leads to an interesting problem: in \cite{Thiel-Dec} we showed that decomposition matrices of $A$ are trivial precisely on an open subset $\msf{DecGen}(A)$ of $\Spec(R)$. In Theorem \ref{decgen_blgen_inclusion_main_theorem} we show:
\[
\msf{DecGen}(A) \subs \msf{BlGen}(A) \;.
\]
The obvious question is: are these two sets equal, and if not, when are they equal? We show in Example \ref{decgen_eq_blgen_counterex} that in general we do not have equality. In Lemma \ref{cellular_alg_blgen_decgen}, on the other hand, we establish a context where we have equality (this includes Brauer algebras and explains why our Weil divisor is given by the non-semisimple parameters). It is an open problem to understand the complement $\msf{BlGen}(A) \setminus \msf{DecGen}(A)$.

\subsection{Semicontinuity of blocks in case of a non-split generic fiber}

In \S\ref{non_split_semicont} we consider the case of a non-split generic fiber. In this case we can no longer identify blocks of specializations with blocks of localizations, and so there is no natural way of comparing block structures among the fibers. However, it still makes sense to compare the \textit{number} of blocks of the fibers, i.e., to consider the map $\Spec(R) \rarr \bbN$, $\fp \mapsto \#\msf{Bl}(A(\fp))$. In case $R$ is normal and $A^K$ splits, this map is lower semicontinuous by the results discussed above. Without splitting of $A^K$, this is no longer true, see Example \ref{ken_example}. The problem is that we consider this map on all of $\Spec(R)$. In Corollary \ref{bln_spec_strat} we construct a setting in which the restriction of $\fp \mapsto \#\msf{Bl}(A(\fp))$ to certain subsets of $\Spec(R)$ is still lower semicontinuous without assuming that the generic fiber splits. From this we obtain a rather nice result, see  Corollary \ref{finite_type_alg_closed_setting}:

\begin{thm}  \label{non-split-intro}
Suppose that $R$ is a finite type algebra over an algebraically closed field. Let $X$ be the set of closed points of $\Spec(R)$. Then the map $X \rarr \bbN$, $\fm \mapsto \#\msf{Bl}(A(\fm))$, is lower semicontinuous. In particular, $X$ admits a stratification according to the number of blocks of fibers of $A$ over $X$.
\end{thm}

\subsection{Remark}
The behavior of blocks under specialization has been studied in several  situations already. All of our results are well-known in modular representation theory of finite groups since the work of R. Brauer and C. Nesbitt \cite{Brauer-Nesbitt}. Our Corollary \ref{normal_split_corollary} and Theorem \ref{dec_block_compat} generalize results by S. Donkin and R. Tange \cite{Donkin-Tange-Brauer} about algebras over Dedekind domains. Our results about lower semicontinuity of the number of blocks generalize a result by P. Gabriel \cite{Gabriel-finite-rep-type-open} to mixed characteristic and non-algebraically closed settings, see also the corresponding result by I. Gordon \cite{Gordon:Representations-of-semisimple-Lie-in-pos}. In general, K. Brown and I. Gordon \cite{BG-Ramification, BG-Ramification2} used Müller's theorem \cite{Mul-Localization-in-non-commu-0} to study blocks under specialization. Theorem \ref{mueller_maintheorem} has been treated in a more special setting by K. Brown and K. Goodearl \cite{Brown-Goodearl-Quantum-Groups}. The codimension one property in Corollary \ref{normal_split_corollary} and Theorem \ref{maintheorem_blex_descr} were proven by C. Bonnafé and R. Rouquier \cite{BR-cellules} in a more special setting. Their work is without doubt one of the main motivations for this paper. Blocks and decomposition matrices of generically semisimple algebras over discrete valuation rings have been studied by M. Geck and G. Pfeiffer \cite{Geck-Pfeiffer}, and more generally by M. Chlouveraki \cite{Chlouveraki:2009aa}. Brauer reciprocity has been studied more generally by M. Geck and R. Rouquier \cite{Geck-Rouquier-Dec}, and by M. Neunhöffer \cite{Neunhoeffer-PhD}. M. Neunhöffer and S. Scherotzke \cite{Neunhoeffer-Scherotzke} have shown generic triviality of $\msf{e}_A^\fp$ over Dedekind domains. 

\subsection*{Acknowledgements}
I would like to thank Cédric Bonnafé for many helpful discussions about this topic, for showing and explaining me the relevant part of the manuscript \cite{BR-cellules} with Raphaël Rouquier, and for providing Example \ref{decgen_eq_blgen_counterex}. The manuscript \cite{BR-cellules} is without doubt one of the main motivations for this paper. I would furthermore like to thank Gwyn Bellamy, Ken Brown, Meinolf Geck, and Gunter Malle for commenting on parts of a preliminary version of this article. Moreover, I thank Ken Brown for providing Example \ref{ken_example}. I was partially supported by the DFG SPP 1489.

\tableofcontents

\section{Base change of blocks} \label{notations}

The basic principle underlying the behavior of blocks in a family of algebras is \textit{base change} of blocks. In this section, we introduce a few basic notions about this principle. Appendix \ref{appendix_blocks_base_change} contains some further material which will later be used in some proofs.

Let us fix some basic notations for block theory. For us, a \word{ring} is always a ring with identity and a \word{module} is always a \textit{left} module unless we explicitly say it is a \textit{right} module. 
Let $A$ be a ring and let $Z$ be its center. If $c$ is a central idempotent of $A$, then $Ac = cA$ is a two-sided ideal of $A$ and at the same time a ring with identity element equal to $c$ (so, not a subring). This yields a bijection between between the set of decompositions of $1 \in A$ into a sum of pairwise orthogonal central idempotents and finite direct sum decompositions of the ring $A$ into non-zero two-sided ideals of $A$ up to permutation of the summands. Such decompositions are in turn in bijection with finite direct product decompositions of the ring $A$ into non-zero rings up to permutation of the factors. Primitive idempotents of $Z$ are also called \word{centrally-primitive} idempotents of $A$. A central idempotent $c$ is centrally-primitive if and only if $Ac$ is an indecomposable ring.
 It is a standard fact—and the starting point of block theory—that if there is a decomposition of $1 = \sum_i c_i$ into pairwise orthogonal centrally-primitive idempotents $c_i$, then this is unique and any central idempotent of $A$ is a sum of a subset of the $c_i$. We then say that $A$ has a \word{block decomposition}, call the centrally-primitive idempotents of $A$ also the \word{block idempotents}, and call the corresponding rings $Ac$ the \word{blocks} of $A$. In this case we prefer to write $\msf{Bl}(A) \dopgleich \msf{Idem}_{cp}(A)$. To avoid pathologies we set $\msf{Bl}(0) \dopgleich \emptyset$ for the zero ring $0$. It is well-known that noetherian rings have block decompositions (the block idempotents are the class sums with respect to the linkage relation of a decomposition of $1 \in A$ into pairwise orthogonal primitive idempotents).

Let $\mscr{C} \dopgleich \lbrace c_i \rbrace_{i \in I}$ be a finite set of pairwise orthogonal central idempotents whose sum is equal to $1 \in A$. Let $B_i \dopgleich Ac_i$. If $V$ is a non-zero $A$-module, then $V = \bigoplus_{i \in I} c_i V$ as $A$-modules and each summand $c_i V$ is a $B_i$-module. In this way we obtain a decomposition $A\tn{-}\msf{Mod} = \bigoplus_{i \in I} B_i\tn{-}\msf{Mod}$ of module categories, which also restricts to a decomposition of the category of finitely generated modules. If a non-zero $A$-module $V$ is under this decomposition obtained from a $B_i$-module, then $V$ is said to \word{belong} to $B_i$. This is equivalent to $c_i V = V$ and $c_j V = 0$ for all $j \neq i$. An indecomposable, and thus any simple, $A$-module clearly belongs to a unique $B_i$. We thus get a decomposition $\msf{Irr} A = \coprod_{i \in I} \msf{Irr} B_i$ of the set of (isomorphism classes of) simple modules. We call the sets $\msf{Irr} B_i$ the \word{$\mscr{C}$-families} of $A$ and denote the set of $\mscr{C}$-families by $\msf{Fam}_{\mscr{C}}(A)$. Note that we have a natural bijection
\begin{equation}
\mscr{C} \overset{\sim}{\longrightarrow} \msf{Fam}_{\mscr{C}}(A) 
\end{equation}
given by $c_i \mapsto \Irr B_i$. In case $\mscr{C}$ is actually a block decomposition, we call the $\mscr{C}$-families simply the \word{families} of $A$ and set $\msf{Fam}(A) \dopgleich \msf{Fam}_{\mscr{C}}(A)$. Recall that any central idempotent of $A$ is a sum of a subset of the block idempotents of $A$. Hence, for general $\mscr{C}$ as above the families are a finer partition of $\msf{Irr}A$ than the $\mscr{C}$-families, i.e., any $\mscr{C}$-family is a union of families.

Now, consider a morphism  $\phi:R \rarr S$ of commutative rings. If $V$ is an $R$-module, we write $V^S \dopgleich \phi^*V \dopgleich S \otimes_R V$ for the scalar extension of $V$ to $S$ and by $\phi_V:V \rarr V^S$ we denote the canonical map $v \mapsto 1 \otimes v$. In most situations we consider, this map will be injective:

\begin{lem} \label{phi_injective_lemma}
In each of the following cases the map $\phi_V:V \rarr V^S$ is injective:
\begin{enum_thm}
\item \label{phi_injective_lemma:proj} $\phi$ is injective and $V$ is $R$-projective.
\item \label{phi_injective_lemma:faithfully_flat} $\phi$ is faithfully flat.
\item \label{phi_injective_lemma:localiz} $\phi$ is the localization morphism for a multiplicatively closed subset $\Sigma \subs R$ and $V$ is $\Sigma$-torsion-free.
\end{enum_thm}
\end{lem}

\begin{proof}
The first case follows from \cite[II, \S5.1, Corollary to Proposition 4]{Bou-Algebra-1-3}, the second follows from \cite[I, \S3.5, Proposition 8(i,iii)]{Bou-Commutative-Algebra-1-7}, and the last case follows from the fact that $\phi$ is flat in conjunction with \cite[I, \S2.2, Proposition 4]{Bou-Commutative-Algebra-1-7}.
\end{proof}

If $A$ is an $R$-algebra, then the $S$-module $A^S$ is naturally an $S$-algebra and the map $\phi_A:A \rarr A^S$ is a ring morphism. Moreover, if $V$ is an $A$-module, then the underlying $S$-module of $A^S \otimes_A V$ is simply $V^S$. 
Our aim is to study the behavior of blocks under the morphism $\phi_A:A \rarr A^S$. Clearly, if $e \in A$ is an  idempotent, also $\phi_A(e) \in A^S$ is an idempotent, and if $e$ is central, so is $\phi_A(e)$ by the elementary fact that  
\begin{equation}
\phi_A(\msf{Z}(A)) \subs \msf{Z}(A^S) \;.
\end{equation}

\begin{defn}
We say that $\phi_A$ is (central) \word{idempotent stable} if $\phi_A(e) \neq 0$ for any non-zero (central) idempotent $e$ of $A$. We say that $\phi_A$ is \word{block bijective} if $\phi_A$ induces a bijection between the centrally-primitive idempotents of $A$ and the centrally-primitive idempotents of $A^S$.
\end{defn}

Note that in case $\phi_A$ is idempotent stable, respectively central idempotent stable, it induces a map between the sets of decompositions of $1 \in A$ and $1 \in A^S$ into pairwise orthogonal idempotents, respectively into pairwise orthogonal central idempotents. The following lemma shows two situations in which $\phi_A$ is idempotent stable (and thus central idempotent stable). We denote by $\msf{Rad}(A)$ the Jacobson radical of $A$.

\begin{lem} \label{idempotent_stable}
If $\msf{Ker}(\phi_A) \subs \msf{Rad}(A)$, then $\phi_A$ is idempotent stable. This holds in the following two cases:
\begin{enum_thm}
\item \label{idempotent_stable:inj} $\phi_A$ is injective (see Lemma \ref{phi_injective_lemma}),
\item \label{idempotent_stable:surj} $\phi$ is surjective, $\msf{Ker}(\phi) \subs \msf{Rad}(R)$, and $A$ is finitely generated as an $R$-module.
\end{enum_thm}
\end{lem}

\begin{proof}
If $e \in A$ is an idempotent contained in $\msf{Rad}(A)$, then by a well-known characterization of the Jacobson radical (see \cite[5.10]{CR-Methods-1}) we conclude that $e^\dagger = 1-e \in A^\times$ is a unit, and since $e^\dagger$ is also an idempotent, we must have $e^\dagger = 1$, implying that $e = 0$. If $\phi_A$ is injective, the condition clearly holds. In the second case we have $
\msf{Ker}(\phi_A) = \msf{Ker}(\phi) A \subs \msf{Rad}(R) A \subs \msf{Rad}(A)$, where the last inclusion follows from \cite[Corollary 5.9]{Lam-First-Course-91}.
\end{proof}

Suppose that $\phi_A$ is idempotent stable and that both $A$ and $A^S$ have block decompositions. Let $\lbrace c_i \rbrace_{i \in I}$ be the block idempotents of $A$ and let $\lbrace c_j' \rbrace_{j \in J}$ be the block idempotents of $A^S$. Since $\phi_A$ is idempotent stable, the set $\msf{Bl}_\phi(A^S) \dopgleich \phi_A(\lbrace c_i \rbrace_{i \in I})$ is a decomposition of $1 \in A^S$ into pairwise orthogonal idempotents. We call the $\phi_A(c_i)$ the \word{$\phi$-blocks} of $A^S$ and call the corresponding families (see above) the \word{$\phi$-families} of $A^S$, denoted $\msf{Fam}_\phi(A^S)$. As explained above, each $\phi$-block $\phi_A(c_i)$ is a sum of a subset of the block idempotents of $A^S$ and the $\phi$-families are coarser than the families in the sense that each $\phi$-family is a union of $A^S$-families. In particular, we have
\begin{equation} \label{get_more_blocks_equation}
\# \msf{Bl}(A) = \# \msf{Bl}_\phi(A^S) \leq \# \msf{Bl}(A^S) \;.
\end{equation}
The following picture illustrates this situation:
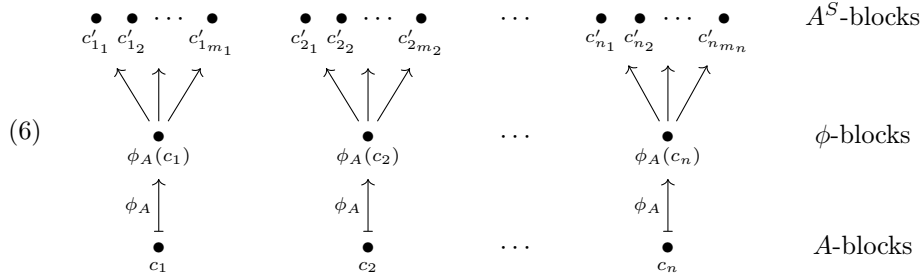
\begin{equation}
\begin{tikzcd}[column sep=small]
\underset{c_{1_1}'}{\bullet}  \underset{c_{1_2}'}{\bullet}  \cdots \underset{c_{1_{m_1}}'}{\bullet} & \underset{c_{2_1}'}{\bullet}  \underset{c_{2_2}'}{\bullet}  \cdots \underset{c_{2_{m_2}}'}{\bullet}  & \cdots & \underset{c_{n_1}'}{\bullet}  \underset{c_{n_2}'}{\bullet}  \cdots \underset{c_{n_{m_n}}'}{\bullet}  & A^S\tn{-blocks} \\
\underset{\phi_A(c_1)}{\bullet} \arrow{u} \arrow[end anchor=230]{u} \arrow[end anchor=310]{u} & \underset{\phi_A(c_2)}{\bullet} \arrow{u} \arrow[end anchor=230]{u} \arrow[end anchor=310]{u}  & \cdots & \underset{\phi_A(c_n)}{\bullet} \arrow{u} \arrow[end anchor=230]{u} \arrow[end anchor=310]{u}  & \phi\tn{-blocks} \\
\underset{c_1}{\bullet} \arrow[mapsto]{u}{\phi_A} & \underset{c_2}{\bullet} \arrow[mapsto]{u}{\phi_A} & \cdots & \underset{c_n}{\bullet} \arrow[mapsto]{u}{\phi_A} & A\tn{-blocks}
\end{tikzcd}
\end{equation}
This paper is about this picture in the special case of specializations of an algebra in prime ideals. Before we begin investigating this, we record the following useful fact.

\begin{lem} \label{block_dec_reflection}
Suppose that $\phi_A:A \rarr A^S$ is central idempotent stable. If $A^S$  has a block decomposition, then $A$ has a block decomposition.
\end{lem}

\begin{proof}
If $A$ does not contain any non-trivial central idempotent, then $A$ is indecomposable and thus has a block decomposition. So, assume that $A$ is not indecomposable and let $c$ be a non-trivial central idempotent. Then $A = Ac \oplus Ac^\dagger$. We can now continue this process to get finer and finer decompositions of $A$ as a ring. Since $\phi_A$ is central idempotent stable, we get decompositions of the same size of $A^S$. As $A^S$ has a block decomposition, this process has to end after finitely many steps. We thus arrive at a ring decomposition of $A$ with finitely many and indecomposable factors, hence, at a block decomposition of $A$.
\end{proof}

\begin{cor} \label{finite_flat_int_block_dec}
A non-zero finite flat algebra over an integral domain has a block decomposition.
\end{cor}

\begin{proof}
Let $R$ be an integral domain with fraction field $K$, let $\phi:R \hookrightarrow K$ be the embedding, and let $A$ be a finite flat $R$-algebra. Since $A$ is $R$-torsion-free, it follows from Lemma \ref{phi_injective_lemma}\ref{phi_injective_lemma:localiz} that $\phi_A$ is injective and so $\phi_A$ is idempotent stable by Lemma \ref{idempotent_stable}\ref{idempotent_stable:inj}. Since $\phi_A^*A = A^K$ is a finite-dimensional algebra over a field, it has a block decomposition. Hence, $A$ has a block decomposition by Lemma \ref{block_dec_reflection}.
\end{proof}

The point of the corollary above is that we do \textit{not} have to assume $R$ to be noetherian—otherwise $A$ is noetherian and we already know it has a block decomposition.

\section{Blocks of localizations} \label{blocks_of_localizations}

Before we consider blocks of specializations, we first take a look at blocks of localizations as these are much easier to control and are still strongly related to blocks of specializations as we will see in the next paragraph. %

\begin{tcolorbox}
Throughout this paragraph, we assume that $A$ is a finite flat algebra over an integral domain $R$ with fraction field $K$. 
\end{tcolorbox}

It follows from Corollary \ref{finite_flat_int_block_dec} that $A$ and any localization $A_\fp$ for $\fp \in \Spec(R)$ has a block decomposition, even if $A$ is not necessarily noetherian. Since the canonical map $\phi_\fp: A_\fp \rarr A^K$ is injective by Lemma \ref{phi_injective_lemma}, we have the notion of $\phi_\fp$-blocks and $\phi_\fp$-families of $A^K$ as defined in \S\ref{notations}. To shorten notations, we call them the \word{$\fp$-blocks} and \word{$\fp$-families}, and write $\msf{Fam}_\fp(A^K)$ for the $\fp$-families. Recall that we have a natural bijection 
\begin{equation}
\msf{Bl}(A_\fp) \simeq \msf{Fam}_\fp(A^K) \;.
\end{equation}

\subsection{Block structure stratification}

There is the following more concrete point of view of $\fp$-blocks. Let $(c_i)_{i \in I}$ be the block idempotents of $A^K$. If $c \in A_\fp$ is any block idempotent, we know from \S\ref{notations} that there is $I' \subs I$ with $c = \sum_{i \in I'} c_i$ in $A^K$. Hence, to any block idempotent of $A_\fp$ we can associate a subset of $I$, and if we take all block idempotents of $A_\fp$ into account, we get a partition $\Bl_A(\fp)$ of the set $I$, from which we can recover the block idempotents of $A_\fp$ by taking sums of the $c_i$ over the members of $\Bl_A(\fp)$. In this way we get a map
\begin{equation}
\Bl_A:\Spec(R) \rarr \msf{Part}(I) 
\end{equation}
to the set of partitions of the set $I$. We denote by
\begin{equation}
\Bl(A) \dopgleich \Im \Bl_A
\end{equation}
the image of this map and call the partitions therein the \word{block structures} of $A$. 

The set $\msf{Part}(I)$ is equipped with the partial order $\leq$ defined by $\mscr{P} \leq \mscr{Q}$ if $\mscr{P}$ is a coarser partition than $\mscr{Q}$, i.e., the members of $\mscr{P}$ are unions of members of $\mscr{Q}$. If $\fq \subs \fp$, then we have an embedding $A_\fp \hookrightarrow A_\fq$ and by the same argumentation as above, the block idempotents of $A_\fp$ are obtained by summing up block idempotents of $A_\fq$, so
\begin{equation}
\fq \subseteq \fp \Longrightarrow \Bl_A(\fp) \leq \Bl_A(\fq) \;.
\end{equation}
 Hence, the map $\Bl_A$ is actually a morphism of posets if we equip $\Spec(R)$ with the partial order $\leq$ defined by $\fp \leq \fq$ if $\fq \subs \fp$ (i.e., $\msf{V}(\fp) \subs \msf{V}(\fq)$). 

For $\mscr{P} \in \msf{Part}(I)$, we call the fiber $\Bl_A^{-1}(\mscr{P}) \subs \Spec(R)$ the \word{$\mscr{P}$-stratum} and we call 
\begin{equation} \label{skeleton_def}
\Bl_A^{-1}(\leq \! \mscr{P}) \dopgleich \bigcup_{\mscr{P}' \leq \mscr{P}} \Bl_A^{-1}(\mscr{P}') = \bigcup_{\substack{\mscr{P}' \leq \mscr{P} \\ \mscr{P}' \in \Bl(A)}} \Bl_A^{-1}(\mscr{P}') 
\end{equation}
the \word{$\mscr{P}$-skeleton}. This is simply the locus of all $\fp \in \Spec(R)$ where the block structure of $A_\fp$ is equal to $\mscr{P}$, respectively coarser than $\mscr{P}$. 
Since 
\begin{equation} \label{stratum_from_skeleton}
\Bl_A^{-1}(\mscr{P}) = \Bl_A^{-1}(\leq \! \mscr{P}) \setminus \bigcup_{\mscr{P}' <\mscr{P}} \Bl_A^{-1}(\leq \! \mscr{P'}) \;,
\end{equation}
 we can recover the strata from the skeleta. We get a finite decomposition
\begin{equation} \label{block_stratification}
\Spec(R) = \coprod_{\mscr{P}} \Bl_A^{-1}(\mscr{P})  
\end{equation}
and we call this the \textit{block structure stratification}. 
Our aim is now to show that this is indeed a stratification, i.e., the strata are \textit{locally closed} subsets of $\Spec(R)$ and the closure of a stratum is contained in its skeleton. The key ingredient in proving this is the following general proposition, which is essentially due to Bonnafé and Rouquier \cite[Proposition D.2.11]{BR-cellules} but is proven here in a more general form.

\begin{prop} \label{br_gen_lemma}
Let $R$ be an integral domain with fraction field $K$, let $A$ be a finite flat $R$-algebra, and let $\mscr{F} \subs A^K$ be a finite set. Then
\[
\msf{Gen}_A(\mscr{F}) \dopgleich \lbrace \fp \in \Spec(R) \mid \mscr{F} \subs A_\fp \rbrace
\]
is a \textnormal{neighborhood} of the generic point of $\Spec(R)$. If $A$ is finitely presented flat, then $\msf{Gen}_A(\mscr{F})$ is an \textnormal{open} subset of $\Spec(R)$, and if moreover $R$ is a Krull domain, the complement $\msf{Ex}_A(\mscr{F})$ of $\msf{Gen}_A(\mscr{F})$ in $\Spec(R)$ is a \textnormal{reduced Weil divisor}, i.e., it is either empty or pure of codimension one with finitely many irreducible components.
\end{prop}

\begin{proof}
Let us first assume that $A$ is actually $R$-free. For an element $\alpha \in K$ we define $I_\alpha \dopgleich \lbrace r \in R \mid r\alpha \in R \rbrace$. This is a non-zero ideal in $R$, and it has the property that $\alpha \in R_{\fp }$ if and only if $I_\alpha \nsubseteq \fp $. To see this, suppose that $\alpha \in R_{\fp }$. Then we can write $\alpha = \frac{r}{x}$ for some $x \in R \setminus \fp $. Hence, $x \alpha = r \in R$ and therefore $x \in I_\alpha$. Since $x \notin \fp $, it follows that $I_\alpha \nsubseteq \fp $. Conversely, if $I_\alpha \nsubseteq \fp $, then there exists $x \in I_\alpha$ with $x \notin \fp $. By definition of $I_\alpha$ we have $x \alpha \gleichdop r \in R$ and since $x \notin \fp $, we can write $\alpha = \frac{r}{x} \in R_{\fp }$. Now, let $(a_1,\ldots,a_n)$ be an $R$-basis of $A$. Then we can write every element $f \in \sF$ as $f = \sum_{i=1}^n \alpha_{f,i} a_i$ with $\alpha_{f,i} \in K$. Let $I$ be the radical of the ideal
\[
\prod_{{f \in \sF, \; i =1,\ldots,n}} I_{\alpha_{f,i}} \unlhd R \;.
\] 
By the properties of the ideals $I_\alpha$ we have the following logical equivalences:
\[
\begin{array}{rcl}
(\sF \subs A_{\fp }) & \Longleftrightarrow & (\alpha_{f,i} \in R_{\fp } \quad \forall f \in \sF, \ i =1,\ldots,n) \\ &\Longleftrightarrow& (I_{\alpha_{f,i}} \not\subs \fp  \quad \forall f  \in \sF, \ i=1,\ldots,n) \\ & \Longleftrightarrow & (I \not\subs \fp ) \;,
\end{array}
\]
the last equivalence following from the fact that $\fp $ is prime. Hence, 
\begin{equation}
\msf{Ex}_A(\sF) = \Spec(R) \setminus \msf{Gen}_A(\sF) = \msf{V}(I) = \bigcup_{f \in \mscr{F}, \; i=1,\ldots,n} \msf{V}(I_{\alpha_{f,i}}) \;,
\end{equation}
implying that $\msf{Gen}_A(\sF)$ is an open subset of $\Spec(R)$.

Next, still assuming that $A$ is $R$-free, suppose that $R$ is a Krull domain. To show that $\msf{Ex}_A(\mscr{F})$ is either empty or pure of codimension $1$ in $\Spec(R)$ with finitely many irreducible components, it suffices to show this for the closed subsets $\msf{V}(I_\alpha) = \msf{V}(\sqrt{I}_\alpha)$. If $\alpha \in R$, then $I_\alpha = R$ and therefore $\msf{V}(I_\alpha) = \emptyset$. So, let $\alpha \notin R$. Let $\msf{V}(I_\alpha) = \bigcup_{\lambda \in \Lambda} \msf{V}(\fq_\lambda)$ be the decomposition into irreducible components. Note that this decomposition is unique and contains every irreducible component of $\msf{V}(I_\alpha)$. The inclusion $\msf{V}(I_\alpha) \sups \msf{V}(\fq_\lambda)$ is equivalent to $I_\alpha \subs \sqrt{I_\alpha} \subs \sqrt{\fq_\lambda} = \fq_\lambda$. Since an irreducible component is a maximal proper closed subset, we see that the $\fq_\lambda$ are the minimal prime ideals of $\Spec(R)$ containing $I_\alpha$. Let $\fq = \fq_\lambda$ for an arbitrary $\lambda \in \Lambda$. We will show that $\msf{ht}(\fq) = 1$. Since $I_\alpha \subs \fq$, we have seen above that $\alpha \notin R_{\fq}$. As $R$ is a Krull domain, also $R_{\fq}$ is a Krull domain by \cite[Theorem 12.1]{Mat-Commutative}. By \cite[VII, \S1.6, Theorem 4]{Bou-Commutative-Algebra-1-7} we have 
\[
R_{\fq} = \bigcap_{\substack{ \fq' \in \Spec(R_\fq) \\ \msf{ht}(\fq') = 1}} (R_\fq)_{\fq'} = \bigcap_{\substack{ \fq' \in \Spec(R) \\ \fq' \subs \fq \\ \msf{ht}(\fq') = 1}} R_{\fq'} \;.
\]
Since $\alpha \notin R_{\fq}$, this shows that there exists $\fq' \in \Spec(R)$ with $\fq' \subs \fq$, $\msf{ht}(\fq') = 1$ and $\alpha \notin R_{\fq'}$. The last property implies $I_\alpha \subs \fq'$ and now the minimality in the choice of $\fq$ implies that $\fq' = \fq$. Hence, $\msf{ht}(\fq) = 1$ and this shows $\msf{V}(I_\alpha)$ is pure of codimension $1$. Since $I_\alpha \neq 0$, there is some $0 \neq r \in I_\alpha$. This element is contained in all the height one prime ideals $\fq_\lambda$. As $R$ is a Krull domain, a non-zero element of $R$ can only be contained in finitely many height one prime ideals (see \cite[4.10.1]{Huneke-Swanson-Integral-Closure}), so $\Lambda$ must be finite.

Now, assume that $R$ is an arbitrary integral domain and that $A$ is finite flat. Then Grothendieck's generic freeness lemma \cite[Lemme 6.9.2]{Grothendieck:EGA-4-2} shows that there exists a non-zero $f \in R$ such that $A_f$ is a free $R_f$-module. Note that $\Spec(R_f)$ can be identified with the distinguished open subset $\msf{D}(f)$ of $\Spec(R)$. We obviously have
\[
\msf{Gen}_{A_f}(\mscr{F}) = \msf{Gen}_A(\mscr{F}) \cap \msf{D}(f) \;.
\]
By the arguments above, $\msf{Gen}_{A_f}(\mscr{F})$ is an open subset of $\msf{D}(f)$, and thus of $\Spec(R)$. This shows that $\msf{Gen}_A(\mscr{F})$ is a neighborhood in $\Spec(R)$.

Next, let $R$ be arbitrary and assume that $A$ is finitely presented flat. It is a standard fact (see \cite[Tag 00NX]{stacks-project}) that the assumptions on $A$ imply that $A$ is already finite locally free, i.e., there exist a family $(f_i)_{i \in I}$ of elements of $R$ such that the standard open affines $\msf{D}(f_i)$ cover $\Spec(R)$ and $A_{f_i}$ is a finitely generated free $R_{f_i}$-module for all $i \in I$. Since $\Spec(R)$ is quasi-compact, see \cite[Proposition 2.5]{GorWed10-Algebraic-geomet}, we can assume that $I$ is finite. Again note that $\Spec(R_{f_i})$ can be identified with $\msf{D}(f_i)$ and that
\begin{equation} \label{pure_codim_proof_1}
\msf{Gen}_{A_{f_i}}(\mscr{F}) = \msf{Gen}_A(\mscr{F}) \cap \msf{D}(f_i) \;.
\end{equation}
By the above, the set $\msf{Gen}_{A_{f_i}}(\mscr{F})$ is open and since the $\msf{D}(f_i)$ cover $\Spec(R)$, it follows that $\msf{Gen}_A(\mscr{F})$ is open. Now, suppose that $R$ is a Krull domain. Similarly as in (\ref{pure_codim_proof_1}) we have
\begin{equation}
\msf{Ex}_{A_{f_i}}(\mscr{F}) = \msf{Ex}_A(\mscr{F}) \cap \msf{D}(f_i) \;.
\end{equation}
Suppose that $\msf{Ex}_A(\mscr{F})$ is not empty and let $Z$ be an irreducible component of $\msf{Ex}_A(\mscr{F})$. There is an $i \in I$ with $Z \cap \msf{D}(f_i) \neq \emptyset$. The map $T \mapsto \ol{T}$ defines a bijection between irreducible closed subsets of $\msf{D}(f_i)$ and irreducible closed subsets of $\Spec(R)$ which meet $\msf{D}(f_i)$, see \cite[\S1.5]{GorWed10-Algebraic-geomet}. This implies that $Z \cap \msf{D}(f_i)$ is an irreducible component of $\msf{Ex}_A(\mscr{F}) \cap \msf{D}(f_i) = \msf{Ex}_{A_{f_i}}(\mscr{F})$. It follows from the above that $Z \cap \msf{D}(f_i)$ is of codimension $1$ in $\msf{D}(f_i)$. Hence, $Z$ is of codimension $1$ in $\Spec(R)$ by \cite[Tag 02I4]{stacks-project}. All irreducible components of $\msf{Ex}_A(\mscr{F})$ are thus of codimension $1$ in $\Spec(R)$. Since each set $\msf{Ex}_{A_{f_i}}(\mscr{F})$ has only finitely many irreducible components and since $I$ is finite, also $\msf{Ex}_A(\mscr{F})$ has only finitely many irreducible components.
\end{proof}

\begin{rem} \label{fp_projective_remark}
We note that $A$ is finitely presented flat if and only if it is finite projective, see \cite[Theorem 4.30]{Lam-Lectures-Modules-Rings-99} or \cite[Tag 058R]{stacks-project}. Hence, we could have equally assumed that $A$ is finite projective in Proposition \ref{br_gen_lemma} but we preferred the seemingly more general notion.
\end{rem}

\begin{tcolorbox}
From now on, we assume that $A$ is finitely presented as an $R$-module.
\end{tcolorbox}

For $\fp \in \Spec(R)$ let us denote by $\mscr{B}_A(\fp) \subs A^K$ the set of block idempotents of $A_\fp$. Clearly, $\mscr{B}_A(\fp)$ and $\Bl_A(\fp)$ are in bijection by taking sums of the $c_i$ over the subsets in $\Bl_A(\fp)$. Note that $\mscr{B}_A(\fp)$ is constant on $\Bl_A^{-1}(\mscr{P})$ for any $\mscr{P}$. We can thus define $\msf{Gen}_A(\mscr{P}) \dopgleich \msf{Gen}_A(\mscr{B}_A(\fp))$ where $\fp \in \Bl_A^{-1}(\mscr{P})$ is arbitrary. 

\begin{thm} \label{block_stratification_thm}
Then $\Bl_A^{-1}(\mscr{P})$ is a closed subset of $\Spec(R)$ for any partition $\mscr{P}$. Hence, each stratum $\Bl_A^{-1}(\mscr{P})$ is \textnormal{open} in $\Bl_A^{-1}(\leq \! \mscr{P})$, thus \textnormal{locally closed} in $\Spec(R)$, and
\begin{equation}
\ol{\Bl_A^{-1}(\mscr{P})} \subs \Bl_A^{-1}(\leq \! \mscr{P}) \;.
\end{equation}
In particular, the decomposition \eqref{block_stratification} is a \textnormal{stratification} of $\Spec(R)$.
\end{thm}

\begin{proof}
First, assume that $\mscr{P}$ is actually a block structure, i.e., $\mscr{P} \in \Bl(A)$. Since $\Spec(R) = \coprod_{\mscr{P}'} \Bl_A^{-1}(\mscr{P}')$, we have 
\[
\Spec(R) \setminus \Bl_A^{-1}(\leq \! \mscr{P}) = \bigcup_{\mscr{P}' \not\leq \mscr{P}} \Bl_A^{-1}(\mscr{P}') \;.
\]
Let $\mscr{P}' \not\leq \mscr{P}$ and $\fp' \in \msf{Gen}_A(\mscr{P}')$. Then $\mscr{P}' \leq \Bl_A(\fp')$. But this implies that $\Bl_A(\fp') \not\leq \mscr{P}$ since otherwise $\mscr{P}' \leq \Bl_A(\fp') \leq \mscr{P}$. Hence, $\msf{Gen}_A(\mscr{P}') \subs \Spec(R) \setminus \Bl_A^{-1}(\leq \! \mscr{P})$. Conversely, we clearly have $\Bl_A^{-1}(\mscr{P}') \subs \msf{Gen}_A(\mscr{P}')$. This shows that
\[
\Spec(R) \setminus \Bl_A^{-1}(\leq \! \mscr{P}) = \bigcup_{\mscr{P}' \not\leq \mscr{P}} \Bl_A^{-1}(\mscr{P}') = \bigcup_{\mscr{P}' \not\leq \mscr{P}} \msf{Gen}_A(\mscr{P}') \;.
\]
This set is open by Proposition \ref{br_gen_lemma}, so 
\begin{equation}
\Bl_A^{-1}(\leq \! \mscr{P}) = \bigcap_{\mscr{P}' \not\leq \mscr{P}} \msf{Ex}_A(\mscr{P}')
\end{equation}
is closed. From \eqref{skeleton_def} we see that $\Bl_A^{-1}(\leq \mscr{P})$ is also closed for an arbitrary partition $\mscr{P}$. Using (\ref{stratum_from_skeleton}), we see that $\Bl_A(\mscr{P})$ is locally closed. Moreover, we have $\Bl_A^{-1}(\mscr{P}) \subs \Bl_A^{-1}(\leq \! \mscr{P})$, so
\[
\ol{\Bl_A^{-1}(\mscr{P})} \subs \ol{\Bl_A^{-1}(\leq \! \mscr{P})} = \Bl_A^{-1}(\leq \! \mscr{P}) \;.
\]
\end{proof}

\begin{rem}
In general it is \textit{not} true that we have equality $\ol{\Bl_A^{-1}(\mscr{P})} = \Bl_A^{-1}(\leq \! \mscr{P})$, so the stratification (\ref{block_stratification}) is in general not a so-called \textit{good} stratification: in the Brauer algebra example in the introduction we have $\mscr{P}' \dopgleich \Bl_A( (3) ) = \lbrace \lbrace 1,2,3 \rbrace, \lbrace 4 \rbrace \rbrace < \lbrace \lbrace 1,2 \rbrace, \lbrace 3 \rbrace, \lbrace 4 \rbrace \rbrace = \Bl_A((2)) \gleichdop \mscr{P}$, so $(3) \in \Bl_A^{-1}(\leq\! \mscr{P})$, but $(3)$ is not contained in $\ol{\Bl_A^{-1}(\mscr{P})} = \msf{V}((2))$. The problem here is that the skeleton $\Bl_A^{-1}(\leq \! \mscr{P})$ has an irreducible component on which the maximal block structure is strictly smaller than the maximal one on the entire skeleton.
\end{rem}

The poset $\Bl(A)$ has a unique maximal element, namely the block structure $\Bl_A(\bullet)$ of $A$ in the generic point $\bullet \dopgleich (0)$ of $\Spec(R)$, i.e., $\Bl_A(\bullet) = \lbrace \lbrace i \rbrace \mid i \in I \rbrace$ is the block structure of the generic fiber $A^K = A_\bullet$. The deviation of block structures from the generic one thus takes place on the closed subset
\begin{equation} \label{blex_def}
\msf{BlEx}(A) \dopgleich \lbrace \fp \in \Spec(R) \mid \Bl_A(\fp) < \Bl_A(\bullet) \rbrace = \bigcup_{\mscr{P} < \Bl_A(\bullet) } \Bl_A^{-1}(\leq \! \mscr{P}) \;.
\end{equation}
We call this set the \word{block structure divisor} of $A$. In fact, since $\msf{BlEx}(A) = \msf{Ex}_A(\mscr{B}_A(\bullet))$, Proposition \ref{br_gen_lemma}  implies:

\begin{cor} \label{blex_is_weil_div}
Suppose that $R$ is a \textnormal{Krull domain}. Then $\msf{BlEx}(A)$ is a reduced Weil divisor.
\end{cor}

The generic block structure lives precisely on the open dense subset
\begin{equation} \label{blgen_def}
\msf{BlGen}(A) \dopgleich \Spec(R) \setminus \msf{BlEx}(A) = \lbrace \fp \in \Spec(R) \mid \Bl_A(\fp) = \Bl_A(\bullet) \rbrace  = \Bl_A^{-1}(\bullet) \;.
\end{equation}

\subsection{Block number stratification} \label{sec_block_num_strat}

From the map $\Bl_A: \Spec(R) \to \msf{Part}(I)$ we obtain the numerical invariant
\begin{equation} \label{block_number_map}
\#\Bl_A : \Spec(R) \to \bbN , \; \fp \mapsto \#\Bl_A(\fp) = \# \Bl(A_\fp) \;.
\end{equation}
This map is again a morphism of posets, so
\begin{equation}
\fq \subseteq \fp \Rightarrow \#\Bl_A(\fp) \leq \#\Bl_A(\fq) \;.
\end{equation}
 For $n \in \bbN$ we have 
\begin{equation}
\#\Bl_A^{-1}(n) = \lbrace \fp \in \Spec(R) \mid \#\Bl(A_\fp)  = n \rbrace
\end{equation}
and we get a decomposition 
\begin{equation} \label{block_number_partition}
\Spec(R) = \coprod_{n \in \bbN} \#\Bl_A^{-1}(n) \;.
\end{equation}
We call this the \word{block number stratification}. This decomposition is of course coarser than the one defined by the fibers of $\Bl_A$. We define
\begin{equation}
\#\Bl_A^{-1}(\leq \! n) \dopgleich \bigcup_{m \leq n} \#\Bl_A^{-1}(m) = \lbrace \fp \in \Spec(R) \mid \# \Bl(A_\fp) \leq n \rbrace \;.
\end{equation}
Since
\begin{equation} \label{bl_leqn_closed}
\#\Bl_A^{-1}(\leq \! n) = \bigcup_{ \#\mscr{P}\leq n } \Bl_A^{-1}(\mscr{P}) \;,
\end{equation}
this set is closed in $\Spec(R)$ by Theorem \ref{block_stratification_thm}. This means that the map $\# \Bl_A : \Spec(R) \rarr \bbN$ is \textnormal{lower semicontinuous}. Hence, $\#\Bl_A^{-1}(n)$ is open in $\#\Bl_A^{-1}(\leq \! n)$, thus locally closed in $\Spec(R)$, and
\begin{equation} \label{bln_closure}
\ol{\#\msf{Bl}_A^{-1}( n)} \subs \#\Bl_A^{-1}(\leq \! n) \;,
\end{equation}
In particular, the partition (\ref{block_number_partition}) is a stratification of $\Spec(R)$. Again, in general it will not be a \textit{good} stratification. Note that
\begin{equation}
\msf{BlEx}(A) = \lbrace \fp \in \Spec(R) \mid \#\Bl_A(\fp) < \#\Bl_A(\bullet) \rbrace = \Bl_A^{-1}(\leq \! \#\Bl_A(\bullet)-1) \;.
\end{equation}
 
\subsection{Block stratification} \label{block_strat}

The poset $\msf{Part}(I)$ is actually a lattice, i.e., it has \textit{meets} $\wedge$ and \textit{joins} $\vee$. The meet $\mscr{P} \wedge \mscr{P}'$ of two partitions is the finest partition of $I$ being coarser than both $\mscr{P}$ and $\mscr{P}'$, and this is obtained by joining members with non-empty intersection. The maximal elements in $\msf{Part}(I)$ not equal to the maximal element itself (the trivial partition) are the partitions $\{ i,j \} \cup (I \setminus \{ i,j \})$ with $i \neq j$. We call these the \textit{atoms} of $\msf{Part}(I)$ and we denote by $\msf{At}(I)$ the set of atoms. This terminology comes from the fact that an arbitrary partition $\mscr{P}$ is the meet of all atoms lying above it:
\[
\mscr{P} = \bigwedge_{ \substack{T \in \msf{At}(I) \\ \mscr{P} \leq T}} T \;.
\]
Because of this property, we say that $\msf{Part}(I)$ is \word{atomic}.

The poset of block structures of $A$ has a similar atomic character. For $i,j \in I$ with $i \neq j$ let us write
\begin{equation}
\msf{Gl}_A(\{ i,j \}) \dopgleich \Bl_A^{-1}(\leq \! \{i,j\} \cup (I \setminus \{ i,j \})) \;.
\end{equation}
This is the locus of all $\fp \in \Spec(R)$ such that the block idempotents $c_i$ and $c_j$ belong to the same block of $A_\fp$, i.e., they are ``glued'' over $\fp$. We thus call this set  a \word{gluing locus}. By Theorem \ref{block_stratification_thm} it is a closed subset of $\Spec(R)$. It is clear that
\begin{equation}
\msf{BlEx}(A) = \bigcup_{i \neq j} \msf{Gl}_A(\{ i,j \}) \;.
\end{equation}
Let $\msf{At}(A)$ be the set of maximal elements of the set of irreducible components of the gluing loci, ordered by inclusion. Then we still have
\begin{equation}
\msf{BlEx}(A) = \bigcup_{Z \in \msf{At}(A)} Z \;.
\end{equation}

\begin{lem}
The $Z \in \msf{At}(A)$ are precisely the irreducible components of $\msf{BlEx}(A)$. 
\end{lem}

\begin{proof}
Let $Y$ be an irreducible component of $\msf{BlEx}(A)$ and let $\xi$ be the generic point of $Y$. Since $\Bl_A(\xi)$ is not the trivial partition, there is $i \neq j$ with $\Bl_A(\xi) \leq \{i,j\} \cup (I \setminus \{i,j \})$. Hence, $Y \subs \msf{Gl}_A(\{i,j\})$. Since $Y$ is a maximal irreducible closed subset of $\msf{BlEx}(A)$ and $\msf{Gl}_A(\{i,j\}) \subs \msf{BlEx}(A)$, it is also a maximal irreducible closed subset of $\msf{Gl}_A(\{i,j\})$, thus equal to an irreducible component $Z$ of $\msf{Gl}_A(\{i,j\})$. It is clear that $Z \in \msf{At}(A)$. Conversely, let $Z \in \msf{At}(A)$. Since $Z \subs \msf{BlEx}(A)$, there is an irreducible component $Y$ of $\msf{BlEx}(A)$ containing $Z$. With the same argumentation as above, there is $Z' \in \msf{At}(A)$ with $Y \subs Z'$. Hence, $Z \subs Y \subs Z'$, and therefore $Z=Y$ by maximality of the elements in $\msf{At}(A)$.
\end{proof}

It now follows that for any $\fp \in \Spec(R)$ we have
\begin{equation}
\Bl_A(\fp) = \bigwedge_{ \substack{\mscr{P} \in \msf{At}(I) \\ \Bl_A(\fp) \leq \mscr{P}}} \mscr{P} = \bigwedge_{ \substack{ Z \in \msf{At}(A)\\ \fp \in Z }} \Bl_A(Z) \;,
\end{equation}
where $\Bl_A(Z)$ denotes the block structure in the generic point of $Z$, i.e., the unique maximal block structure on $Z$. Hence, any block structure of $A$ is a meet of \word{atomic block structures} $\msf{Bl}_A(Z)$ for $Z \in \msf{At}(A)$. Recall from Corollary \ref{blex_is_weil_div} that if $R$ is a Krull domain, the $Z \in \msf{At}(A)$ are all of codimension one in $\Spec(R)$.

Following this observation, we introduce a refined stratification of $\Spec(R)$. For a subset $\mscr{Z} \subs \msf{At}(I)$ we define
\begin{equation}
\msf{Bl}_A^{-1}(\leq \! \mscr{Z}) \dopgleich \bigcap_{Z \in \mscr{Z}} Z
\end{equation}
and 
\begin{equation}
\msf{Bl}_A^{-1}(\mscr{Z}) \dopgleich \bigcap_{Z \in \mscr{Z}} Z \setminus \bigcup_{Z \in \msf{At}(A) \setminus \mscr{Z}} Z \;.
\end{equation}
It is clear that $\msf{Bl}_A^{-1}(\leq \! \mscr{Z})$ is closed in $\Spec(R)$, that $\msf{Bl}_A^{-1}(\mscr{Z}) $ is locally closed in $\Spec(R)$ and that the block structure on $\msf{Bl}_A^{-1}(\mscr{Z})$ is in any point equal to $\wedge_{Z \in \mscr{Z}} \Bl_A(Z)$. Note that in this notation $\Bl_A^{-1}(\leq \! \emptyset) = \Spec(R)$ and 
\begin{equation}
\msf{Bl}_A^{-1}(\emptyset) = \Spec(R) \setminus \bigcup_{Z \in \msf{At}(A)} Z = \msf{BlGen}(A) \;.
\end{equation}
Clearly,
\begin{equation}
\overline{\msf{Bl}_A^{-1}(\mscr{Z})} \subseteq \bigcap_{Z \in \mscr{Z}} Z = \Bl_A^{-1}(\leq \!\mscr{Z}) \;,
\end{equation}
so we obtain a stratification
\begin{equation}
\Spec(R) = \coprod_{\mscr{Z} \subs \msf{At}(A)} \Bl_A^{-1}(\mscr{Z}) 
\end{equation}
refining the block structure stratification \eqref{block_stratification}. We call this the \word{block stratification} of $A$.

\section{Blocks of specializations} \label{blocks_of_specializations}

We now turn to our actual problem, namely blocks of specializations of $A$. Compared to blocks of localizations there is in general no possibility to compare the actual block structures of specializations.  
However, there is a rather general setting where blocks of specializations are naturally identified with blocks of localizations, namely when $R$ is normal and $A^K$ splits. In this case we can compare the actual block structures of specializations and all results from the preceding paragraph are actually also results about blocks of specializations. %
For the proof we need the following general result.

\begin{thm} \label{faithfully_flat_ext_block_bij}
Let $\phi: R \hookrightarrow S$ be a faithfully flat morphism of integral domains and let $A$ be a finite flat $R$-algebra. Let $K$ and $L$ be the fraction field of $R$ and $S$, respectively. If $\#\msf{Bl}(A^K) = \#\msf{Bl}(A^L)$, then the morphism $\phi_A:A \rarr A^S$ is block bijective.
\end{thm}

\begin{proof}
Recall from Corollary \ref{finite_flat_int_block_dec} that both $A$ and $A^S$ have block decompositions. The map $\phi_A: A \rarr A^S$ is injective by Lemma \ref{phi_injective_lemma}\ref{phi_injective_lemma:faithfully_flat} since $\phi$ is faithfully flat. Hence, $\phi_A$ is idempotent stable by Lemma \ref{idempotent_stable}\ref{idempotent_stable:inj} and therefore $\# \msf{Bl}(A) \leq \#\msf{Bl}(A^S)$ by (\ref{get_more_blocks_equation}). We thus have to show that $\# \msf{Bl}(A) \geq  \#\msf{Bl}(A^S)$. We split the proof of this fact into several steps. 

The case $R=K$ and $S=L$ holds by assumption. Assume that $R=K$ and that $S$ is general as in the theorem. Since $A$ is $R$-flat, the extension $A^S$ is $S$-flat and thus $S$-torsionfree. Hence, the map $A^S \rarr A^L$ is injective by Lemma \ref{phi_injective_lemma}\ref{phi_injective_lemma:localiz}. In particular, it is idempotent stable by Lemma \ref{idempotent_stable}\ref{idempotent_stable:inj} and so $\#\msf{Bl}(A^S) \leq \#\msf{Bl}(A^L)$ by (\ref{get_more_blocks_equation}). In total, we have $\#\msf{Bl}(A) \leq \# \msf{Bl}(A^S) \leq \# \msf{Bl}(A^L) = \#\msf{Bl}(A^K) = \#\msf{Bl}(A)$. Hence, $\# \msf{Bl}(A) = \# \msf{Bl}(A^S)$.

Finally, let both $R$ and $S$ be general as in the theorem. Let $\Sigma \dopgleich R \setminus \lbrace 0 \rbrace$ and $\Omega \dopgleich S \setminus \lbrace 0 \rbrace$. Then $K = \Sigma^{-1}R$ and $L = \Omega^{-1}S$. Set $T \dopgleich \Sigma^{-1}S$. Since $R$ and $S$ are integral domains, we can naturally view all rings as subrings of $L$ and so we get the two commutative diagrams
\begin{equation}
\begin{tikzcd}[column sep=small, row sep=small]
& L & & & A^L \\
& T \arrow[hookrightarrow]{u} & & & A^T \arrow{u} \\
K \arrow[hookrightarrow]{ur} \arrow[hookrightarrow, bend left]{uur} & & S \arrow[hookrightarrow]{ul} \arrow[hookrightarrow, bend right]{uul} & A^K \arrow[bend left]{uur} \arrow{ur} & & A^S \arrow{ul} \arrow[bend right]{uul} \\
& R \arrow[hookrightarrow]{ur} \arrow[hookrightarrow]{ul} & & & A \arrow{ur} \arrow{ul}
\end{tikzcd}
\end{equation} 
the right one being induced by the left one. All morphisms in the left diagram are clearly injective. We claim the same holds for the right diagram. We have noted at the beginning that the map $A \rarr A^S$ is injective. Since $A$ is $R$-flat, it is $R$-torsionfree and so the map $A \rarr A^K$ is injective by Lemma \ref{phi_injective_lemma}\ref{phi_injective_lemma:localiz}. We have argued above already that the map $A^S \rarr A^L$ is injective. Since $S \hookrightarrow T$ is a localization map, the induced scalar extension functor is exact so that $A^T$ is a flat $T$-module. In particular, $A^T$ is $T$-torsionfree and so $A^T \rarr A^L$ is injective by \ref{phi_injective_lemma}\ref{phi_injective_lemma:localiz}. The map $A^K \rarr A^L$ is injective by Lemma \ref{phi_injective_lemma}\ref{phi_injective_lemma:proj}. Due to the commutativity of the diagram, the remaining maps must be injective, too. We can thus view all scalar extensions of $A$ naturally as subsets of $A^L$. We claim that 
\begin{equation} \label{a_eq_ak_cap_as}
A = A^K \cap A^S
\end{equation}
as subsets of $A^L$. Because of the commutative diagram above, this intersection already takes place in $A^T$. Consider $A^K$ as an $R$-module now. We have a natural identification
\[
\phi^*(A^K) = S \otimes_R A^K = S \otimes_R (\Sigma^{-1} A) = (\Sigma^{-1} S) \otimes_R A = T \otimes_R A = A^T 
\]
as $S$-modules by \cite[II, \S2.7, Proposition 18]{Bou-Commutative-Algebra-1-7}. Note that the map $A^K \rarr A^T$ in the diagram above is the map $\phi_{A^K}$, when considering $A^K$ as an $R$-module. The $R$-submodule $A$ of $A^K$ is now identified with $\phi_{A^K}(A)$ and the $S$-submodule of $A^T$ generated by $A \subs A^T$ is identified with $A^S$. Since $\phi$ is faithfully flat, it follows from \cite[I, \S3.5, Proposition 10(ii)]{Bou-Commutative-Algebra-1-7} applied to the $R$-module $A^K$ and the submodule $A$ that
\[
A = A^K \cap A^S
\]
inside $A^T$. Let $(c_i)_{i \in I}$ be the block idempotents of $A^S$ and let $(d_j)_{j \in J}$ be the block idempotents of $A^K$. By assumption the morphism $A^K \rarr A^L$ is block bijective, which means that $(d_j)_{j \in J}$ are the block idempotents of $A^L$. Since $A^S \rarr A^L$ is idempotent stable, there exists by the arguments preceding  (\ref{get_more_blocks_equation}) a partition $(J_i)_{i \in I}$ of $J$ such that the non-zero central idempotent $c_i$ can in $A^L$ be written as $c_i = \sum_{j \in J_i} d_j$. But this shows that $c_i \in A^K \cap A^S$, hence $c_i \in A$ and so $(c_i)_{i \in I}$ gives a decomposition of $1 \in A$ into pairwise orthogonal centrally-primitive idempotents of $A$ by  (\ref{a_eq_ak_cap_as}). Hence, $\# \msf{Bl}(A) = \# \msf{Bl}(A^S)$.
\end{proof}

To formulate the next proposition more generally, we use the property \textit{block-split} introduced in Definition \ref{block_split_def} but note that the reader might just simply replace it by the more special property \textit{split}. Moreover, we recall that a local integral domain $R$ is called \word{unibranch} if its henselization $R^h$ is again an \textit{integral} (local) domain. This is equivalent to the normalization of $R$ being again \textit{local} (see \cite[IX, Corollaire 1]{Raynaud:Henselian}). This clearly holds if $R$ is already \textit{normal}. Examples of non-normal unibranch rings are the local rings in ordinary cusp singularities of curves. 

\begin{prop} \label{unibranched_block_bijective}
Let $R$ be an integral domain and let $A$ be a finite flat $R$-algebra with \textnormal{block-split} generic fiber $A^K$ (e.g., if $A^K$ splits). Let $\fp \in \Spec(R)$ and suppose that $R_\fp$ is unibranch (e.g., if $R_\fp$ is normal). Then the quotient morphism $A_\fp \twoheadrightarrow A(\fp)$ is block bijective.
\end{prop}

\begin{proof}
By assumption, $R_\fp$ and its henselization $R_\fp^h$ are integral domains. Since $A$ is $R$-flat, it follows that $A_\fp = R_\fp \otimes_R A$ is $R_\fp$-flat and that $A_\fp^h \dopgleich R_\fp^h \otimes_{R_\fp} A_\fp$ is $R_\fp^h$-flat. Hence, both $A_\fp$ and $A_\fp^h$ have block decompositions by Lemma \ref{finite_flat_int_block_dec}. Let $\fp_\fp^h$ be the maximal ideal of $R_\fp^h$. The henselization morphism $R_\fp \rarr R_\fp^h$ is local and faithfully flat by \cite[Théorème 18.6.6(iii)]{Gro67-Elements-de-geom}. We now have a commutative diagram
\[
\begin{tikzcd}
A_\fp \arrow[hookrightarrow]{r} \arrow[twoheadrightarrow,swap]{d} & A_\fp^h \arrow[twoheadrightarrow]{d} \\
A(\fp) = A_\fp/\fp_\fp A_\fp \arrow[hookrightarrow]{r} & A_\fp^h/\fp_\fp^h A_\fp^h
\end{tikzcd}
\]
of idempotent stable morphisms. We know from Lemma \ref{semiperfect_std_settings}\ref{semiperfect_std_settings:hensel} and Lemma \ref{idempotent_surjective_block_bijective} that $A_\fp^h \rarr A_\fp^h/\fp_\fp^h A_\fp^h$ is block bijective. Since $A$ has block-split generic fiber and $R_\fp \rarr R_\fp^h$ is a faithfully flat morphism of integral domains, we can use Theorem \ref{faithfully_flat_ext_block_bij} to deduce that $A_\fp \rarr A_\fp^h$ is block bijective. In \cite[Théorème 18.6.6(iii)]{Gro67-Elements-de-geom} it is proven that $R_\fp/\fp_\fp  \simeq R_\fp^h/\fp_\fp^h $. Hence, the map $A_\fp/\fp_\fp A_\fp \rarr A_\fp^h/\fp_\fp^h A_\fp^h$ is an isomorphism and so in particular block bijective. We thus have
\[
\#\msf{Bl}(A_\fp^h) = \# \msf{Bl}(A_\fp) \leq \# \msf{Bl}(A(\fp))  = \#\msf{Bl}(A_\fp^h/\fp_\fp^h A_\fp^h) = \# \msf{Bl}(A_\fp^h) 
\]
by equation (\ref{get_more_blocks_equation}). Hence, $\# \msf{Bl}(A_\fp) = \# \msf{Bl}(A(\fp))$, so $A_\fp \twoheadrightarrow A(\fp)$ is block bijective.
\end{proof}

\begin{cor} \label{normal_split_corollary}
Suppose that $R$ is normal and $A^K$ splits. Then $A_\fp \twoheadrightarrow A(\fp)$ is block bijective for all $\fp \in \Spec(R)$. Hence, all results from \S\ref{blocks_of_localizations} apply also to blocks of specializations of $A$.
\end{cor}

\section{Blocks via central characters} \label{blocks_via_central_characters}

In this section we discuss an approach to explicitly compute the block structure of $A$ in any point $\fp \in \Spec(R)$, and so to compute the whole block stratification. This is based on the knowledge of the central characters of the generic fiber of $A$. Parts of the arguments presented here are due to Bonnafé and Rouquier \cite[Appendix D]{BR-cellules}.

\subsection{Müller's theorem} \label{sec_mueller}

The central ingredient to establish a relationship between blocks and central characters is the general Lemma \ref{mueller_theorem_idempotents} below, which is usually referred to as \word{Müller's theorem}. We were not able to find a proof of it in this generality in the literature, so we include a proof here but note that this is known. The main ingredient is an even more general result by B. Müller \cite{Mul-Localization-in-non-commu-0} about the fibration of cliques of prime ideals in a noetherian ring over its center, see Lemma \ref{clique_reduction_bijection}. We will recall only a few basic definitions from the excellent exposition in \cite[\S12]{Goodearl-Warfield} and refer to \textit{loc. cit.} for more details.

\begin{tcolorbox}
Throughout this paragraph, we assume that $A$ is a noetherian ring.
\end{tcolorbox}

If $\fp, \fq$ are prime ideals of $A$, we say that there is a \word{link} from $\fp$ to $\fq$, written $\fp \leadsto \fq$, if there is an ideal $\fa$ of $A$ such that $\fp \cap \fq \supsetneq \fa \sups \fp \fq$ and $(\fp \cap \fq)/\fa$ is non-zero and torsion-free both as a left $(A/\fp)$-module and as a right $(A/\fq)$-module. The bimodule $(\fp \cap \fq)/\fq$ is then called a \word{linking bimodule} between $\fq$ and $\fp$. The equivalence classes of the equivalence relation on $\Spec(A)$ generated by $\leadsto$ are called the \word{cliques} of $A$. We write $\msf{Clq}(A)$ for the set of cliques of $A$ and $\msf{Clq}(\fp)$ for the unique clique of $A$ containing $\fp$. For the proof of Lemma \ref{mueller_theorem_idempotents} we will need a few preparatory lemmas. \\

We call the supremum of lengths of chains of prime ideals in a $A$ the \word{classical Krull dimension} of $A$. The following lemma is standard.

\begin{lem} \label{cliques_zero_dim}
Suppose that $A$ is noetherian and of classical Krull dimension zero. Then there is a canonical bijection
\begin{equation}
\begin{array}{rcl} 
\msf{Bl}(A) & \overset{\sim}{\longrightarrow} & \msf{Clq}(A) \\
c & \longmapsto & X_c \dopgleich \lbrace \fm \in \Max(A) \mid c^\dagger \in \fm \rbrace \;,
\end{array}
\end{equation}
where $c^\dagger = 1-c$. If moreover $A$ is \textnormal{commutative}, then the cliques are singletons, i.e., there is a unique $\fm_c \in \Max(A)$ with $c^\dagger \in \fm_c$. Hence, in this case we have $\msf{Bl}(A) \simeq \msf{Max}(A) \simeq \Spec(A)$.
\end{lem}

\begin{proof}
The first assertion is proven in \cite[Corollary 12.13]{Goodearl-Warfield}. In a commutative noetherian ring the cliques are singletons (see \cite[Exercise 12F]{Goodearl-Warfield}), and this immediately implies the second assertion.
\end{proof}

\begin{lem} \label{annihilator_of_torsionfree_mod_p}
Let $\fp $ be a prime ideal of a noetherian ring $A$ and let $V$ be a non-zero $A$-module with $\fp  \subs \msf{Ann}(V)$. If $V$ is torsion-free as an $(A/\fp )$-module, then $\fp  = \msf{Ann}(V)$.
\end{lem}

\begin{proof}
Suppose that $\fp  \subsetneq \msf{Ann}(V)$. Then $\msf{Ann}(V)/\fp $ is a non-zero ideal of the noetherian prime ring $A/\fp $ and thus contains a regular element $\ol{x}$ by \cite[Corollary 2.3.11]{Jat86-Localization-in-}. But then $\ol{x}V = 0$, contradicting the assumption that $V$ is a torsion-free $(A/\fp )$-module.
\end{proof}

\begin{lem} \label{cliques_elementary_notes}
The following holds:
\begin{enum_thm}
\item \label{cliques_elementary_notes:lift} If $\fp $ and $\fq$ are prime ideals of $A$ and if $\fb$ is an ideal of $A$ with $\fb \subs \fp  \cap \fq$ such that $\fp /\fb \leadsto \fq/\fb$ in $A/\fb$, then $\fp  \leadsto \fq$ in $A$.

\item \label{cliques_elementary_notes:reduction} Let $\fp $ and $\fq$ be two prime ideals of $A$ with $\fp  \leadsto \fq$ and let $\fb$ be an ideal of $A$. If there exists a linking ideal $\fa $ from $\fp $ to $\fq$ with $\fb \subs \fa $, then $\fp /\fb \leadsto \fq/\fb$ in $A/\fb$.

\end{enum_thm}
\end{lem}

\begin{proof} \hfill

\begin{enum_proof}

\item We can write a linking ideal from $\fp /\fb$ to $\fq/\fb$ as $\fa /\fb$ for an ideal $\fa $ containing $\fb$. By definition, we have
\[
(\fp  \cap \fq)/\fb = (\fp /\fb) \cap (\fq/\fb) \supsetneq \fa /\fb \sups (\fp /\fb) \cdot (\fq/\fb) = (\fp \fq)/\fb \;,
\]
implying that $\fp  \cap \fq \supsetneq \fa  \sups \fp \fq$. Moreover, we have
\[
\left( ( \fp  \cap \fq)/\fb \right) / \left( \fa/\fb  \right) \cong (\fp  \cap \fq)/\fa 
\]
as $(A/\fb)$-bimodules. By definition, $(\fp  \cap \fq)/\fa $ is torsionfree as a left module over the ring
\[
(A/\fb)/(\fp /\fb) \cong A/\fp \;.
\]
Similarly, it follows that $(\fp  \cap \fq)/\fa $ is torsionfree as a right module over the ring $A/\fq$. Hence, $\fa $ is a linking ideal from $\fp $ to $\fq$.

\smallskip

\item We have
\[
\fp /\fb \cap \fq/\fb = (\fp  \cap \fq)/\fb \supsetneq \fa /\fb \sups (\fp \fq + \fb)/\fb = (\fp /\fb) \cdot (\fq/\fb) \;.
\]
Since 
\[
\left( (\fp  \cap \fq)/\fb \right) / \left( \fa /\fb \right) \cong (\fp  \cap \fq)/\fa \;, \quad (A/\fb)/(\fp /\fb) \cong A/\fp \;, \quad (A/\fb)/(\fq/\fb) \cong A/\fq \;,
\]
it follows that $\fa /\fb$ is a linking ideal from $\fp /\fb$ to $\fq/\fb$. \qedhere
\end{enum_proof}
\end{proof}

\begin{lem}\label{cliques_central_reduction}
Let $\fp $ and $\fq$ be distinct prime ideals of a noetherian ring $A$ with $\fp  \leadsto \fq$. If $\fz$ is a centrally generated ideal of $A$ with $\fz \subs \fp $ or $\fz \subs \fq$, then $\fz \subs \fp  \cap \fq$ and $\fp /\fz \leadsto \fq/\fz$ in $A/\fz$.
\end{lem}

\begin{proof}
This is proven in \cite{Mul85-Affine-Noetheria} but we also give a proof here for the sake of completeness. First note that since $\fz$ is centrally generated and $\fp  \leadsto \fq$, it follows from \cite[Lemma 12.15]{Goodearl-Warfield} that already $\fz \subs \fp  \cap \fq$. Let $\fa $ be a linking ideal from $\fp $ to $\fq$. We claim that $\fz$ is contained in $\fa $. To show this, suppose that $\fz$ is not contained in $\fa $. Then $(\fa  + \fz)/\fa $ is a non-zero submodule of $(\fp  \cap \fq)/\fa $ which is torsionfree as a left $(A/\fp )$-module and as a right $(A/\fq)$-module. In conjunction with the fact that $\fz$ is centrally generated it now follows from Lemma \ref{annihilator_of_torsionfree_mod_p} that
\[
\fp  = \msf{Ann}( _A( (\fa  + \fz)/\fa  ) ) = \msf{Ann}( ((\fa  + \fz)/\fa  )_A) = \fq,
\]
contradicting the assumption $\fp  \neq \fq$. Hence, we must have $\fz \subs \fa $ and it thus follows from Lemma \ref{cliques_elementary_notes}\ref{cliques_elementary_notes:reduction} that $\fp /\fz \leadsto \fq/\fz$. 
\end{proof}

\begin{lem} \label{clique_reduction_bijection}
Let $\fz$ be a centrally generated ideal of a noetherian ring $A$. Let $\fp $ be a prime ideal of $A$ with $\fz \subs \fp $. Then all prime ideals in $\msf{Clq}(\fp )$ contain $\fz$ and the map
\[
\begin{array}{rcl}
\msf{Clq}(\fp )  & \longmapsto & \msf{Clq}(\fp /\fz) \\
\fq & \longmapsto & \fq/\fz
\end{array}
\]
 is a bijection between a clique of $A$ and a clique of $A/\fz$. %
\end{lem}

\begin{proof}
It follows immediately from \cite[Lemma 12.15]{Goodearl-Warfield} that all prime ideals in $\msf{Clq}(\fp )$ contain $\fz$. If $\fq  \in \msf{Clq}(\fp )$, then there exists a chain $\fp  = \fp _0, \fp _1, \ldots, \fp _{r-1}, \fp _r = \fq$ of prime ideals of $A$ with $\fp _i \leadsto \fp _{i+1}$ or $\fp _{i+1} \leadsto \fp _i$ for all indices $i$. An inductive application of Lemma \ref{cliques_central_reduction} shows now that $\fp _i/\fz \leadsto \fp _{i+1}/\fz$ or $\fp _{i+1}/\fz \leadsto \fp _{i}/\fz$ for all $i$. Hence, $\fp /\fz$ and $\fq/\fz$ lie in the same clique of $A/\fz$ so that the map $\msf{Clq}(\fp ) \rarr \msf{Clq}(\fp /\fz)$ is well-defined. On the other hand, similar arguments and Lemma \ref{cliques_elementary_notes}\ref{cliques_elementary_notes:lift} show that if $\fq/\fz \in \msf{Clq}(\fp /\fz)$, then also $\fq \in \msf{Clq}(\fp )$, so that we also have a well-defined map $\msf{Clq}(\fp /\fz) \rarr \msf{Clq}(\fp )$. It is evident that both maps defined are pairwise inverse thus proving the first assertion. The second assertion is now obvious. 
\end{proof}

\begin{lem}[B. Müller] \label{mueller_theorem_idempotents}
Let $A$ be a ring with center $Z$ such that $Z$ is noetherian and $A$ is a finite $Z$-module. If $\fz$ is a centrally generated ideal of $A$ such that $A/\fz A$ is of classical Krull dimension zero, then the inclusion $(Z+\fz)/\fz \hookrightarrow A/\fz A$ is block bijective. In other words, the block idempotents of $A/\fz A$ are already contained in the central subalgebra $(Z+\fz)/\fz$.
\end{lem}

\begin{proof}
Let $\ol{A} \dopgleich A/\fz$ and let $\ol{Z} \dopgleich (Z+\fz)/\fz$. Then $\ol{A}$ is a finitely generated $\ol{Z}$-module since $A$ is a finitely generated $Z$-module. Hence, $\ol{Z} \subs \ol{A}$ is a finite centralizing extension and now it follows from going up in finite centralizing extensions \cite[Theorem 10.2.9]{McR-NN-rings} that the classical Krull dimension of $\ol{Z}$ is equal to that of $\ol{A}$, which is zero by assumption. Hence, by Lemma \ref{cliques_zero_dim} we have $\msf{Bl}(\ol{Z}) \simeq \msf{Clq}(\ol{Z})$ and $\msf{Bl}(\ol{A}) \simeq \msf{Clq}(\ol{A})$. Since $\#\msf{Bl}(\ol{Z}) \leq \#\msf{Bl}(\ol{A})$, the claim is thus equivalent to the claim that over each clique of $\ol{Z}$, there is just one clique of $\ol{A}$. So, let $X,Y \in \msf{Clq}(\ol{A})$ be two cliques. We pick $\fM/\fz \in X$ and $\fN/\fz \in Y$ with $\fM,\fN$ maximal ideals of $A$. Assume that $X$ and $Y$ lie over the same clique of $\ol{Z}$. Since $\ol{Z}$ is commutative, we know from Lemma \ref{cliques_zero_dim} that all cliques are singletons and so the assumption implies that $\fM/\fz$ and $\fN/\fz$ lie over the same maximal ideal of $\ol{Z}$, i.e.,
\[
(\fM/\fz) \cap \left((Z + \fz)/\fz\right)  = (\fN/\fz) \cap \left((Z + \fz)/\fz\right) \;,
\]
hence
\[
\fM \cap (Z+\fz) = \fN \cap (Z+\fz) \;.
\]
Since $Z \subs Z+\fz$, we thus get 
\[
\fM \cap Z = \fM \cap Z \cap (Z+\fz) = \fN \cap Z \cap (Z+\fz) = \fN \cap Z \;.
\]
Now, Müller's theorem \cite[Theorem 13.10]{Goodearl-Warfield} implies that $\fM$ and $\fN$ lie in the same clique of $A$. An application of Lemma \ref{clique_reduction_bijection} thus implies that $\fM/\fz$ and $\fN/\fz$ lie in the same clique of $A/\fz$, so $X=Y$. 
\end{proof}

\subsection{Blocks as fibers of a morphism} \hfill

\begin{tcolorbox}
We assume that $A$ is a finite flat algebra over a \textit{noetherian} integral domain~$R$.
\end{tcolorbox}

By Lemma \ref{ff_center_ext} the morphism 
\begin{equation}
\Upsilon: \Spec(Z) \rarr \Spec(R) \;,
\end{equation}
induced by the canonical morphism from $R$ to the center $Z$ of $A$ is finite, closed, and surjective. The center $Z$ of $A$ is naturally an $R$-algebra and so we can consider its fibers \begin{equation}
Z(\fp) = \msf{k}(\fp) \otimes_R Z/\fp Z = Z_\fp/\fp_\fp Z_\fp
\end{equation}
in prime ideals $\fp$ of $R$. On the other hand, the image of $Z_\fp = \msf{Z}(A_\fp)$ under the canonical (surjective) morphism $A_\fp \twoheadrightarrow A(\fp)$ yields a central subalgebra
\begin{equation}
\msf{Z}_\fp(A) \dopgleich (Z_\fp + \fp_\fp A_\fp)/\fp_\fp A_\fp 
\end{equation}
of $A(\fp)$. In general this subalgebra is \textit{not} equal to the center of $A(\fp)$ itself. We have a surjective morphism 
\begin{equation} \label{center_reduction_map}
\varphi_\fp : Z(\fp) \twoheadrightarrow \msf{Z}_\fp(A) 
\end{equation}
of finite-dimensional $\msf{k}(\fp)$-algebras. This morphism is in general \textit{not} injective—it is if and only if  $\fp_\fp A_\fp \cap Z_\fp = \fp_\fp Z_\fp$. Nonetheless, we have the following result.

\begin{lem} \label{center_reduction_block_bijective}
The map $\varphi_\fp:Z(\fp) \rarr \msf{Z}_\fp(A)$ in (\ref{center_reduction_map}) is block bijective.
\end{lem}

\begin{proof}
Since $\varphi_\fp$ is surjective, the induced map $^a\varphi_\fp: \Spec(\msf{Z}_\fp(A)) \rarr \Spec(Z(\fp))$ is injective, so $\#\Bl(\msf{Z}_\fp(A)) \leq \#\Bl(Z(\fp))$ by Lemma \ref{cliques_zero_dim}. Now we just need to show that $\varphi_\fp$ does not map any no non-trivial idempotent to zero. Since $R_\fp$ is noetherian, also $A_\fp$ is noetherian. The Artin–Rees lemma \cite[Theorem 8.5]{Mat-Commutative} applied to the $R_\fp$-module $A_\fp$, the submodule $Z_\fp$ of $A_\fp$, and the ideal $\fp_\fp$ of $R_\fp$ shows that there is an integer $k \in \bbN_{>0}$ such that for any $n > k$ we have
\[
\fp_\fp^n A_\fp \cap Z_\fp = \fp_\fp^{n-k}(( \fp_\fp^k A_\fp) \cap Z_\fp) \;.
\]
In particular, there is $n \in \bbN_{>0}$ such that $\fp_\fp^n A_\fp \cap Z_\fp \subs \fp_\fp Z_\fp$. Now, let $\ol{e} \in Z(\fp) = Z_\fp/\fp_\fp Z_\fp$ be an idempotent with $\varphi_\fp(\ol{e}) = 0$. By assumption, $\ol{e} \in \Ker(\varphi_\fp) = (\fp_\fp A_\fp \cap Z_\fp)/\fp_\fp Z_\fp $. Hence, if $e \in Z_\fp$ is a representative of $\ol{e}$, we have $e \in \fp_\fp A_\fp \cap Z_\fp$. We have $e^n \in \fp_\fp^n A_\fp \cap Z_\fp \subs \fp_\fp Z_\fp$, so already $\ol{e} = 0$.
\end{proof}

\begin{thm} \label{mueller_maintheorem}
For any $\fp \in \Spec(R)$ there are canonical bijections
\begin{equation} \label{mueller_bijection}
\msf{Bl}(A(\fp)) \simeq \msf{Bl}(\msf{Z}_\fp(A)) \simeq \msf{Bl}(Z(\fp)) \simeq \Upsilon^{-1}(\fp) \;.
\end{equation}
The first bijection $\msf{Bl}(A(\fp)) \simeq \msf{Bl}(\msf{Z}_\fp(A))$ is induced by the embedding $\msf{Z}_\fp(A) \hookrightarrow A(\fp)$. In other words, all block idempotents of $A(\fp)$ are already contained in the central subalgebra $\msf{Z}_\fp(A)$ of $A(\fp)$. The second bijection is the bijection from Lemma \ref{center_reduction_block_bijective}. The last bijection $\msf{Bl}(Z(\fp)) \simeq \Upsilon^{-1}(\fp)$ maps a block idempotent $c$ of $Z(\fp)$ to the (by the theorem unique) maximal ideal $\fm_c$ of $Z$ lying above $\fp$ such that $c^\dagger \in (\fm_c+\fp_\fp Z_\fp)/\fp_\fp Z_\fp$, where $c^\dagger = 1-c$.
\end{thm}

\begin{proof}
The first bijection follows directly from Lemma \ref{mueller_theorem_idempotents} applied to $A_\fp$ and the centrally generated ideal $\fz \dopgleich \fp_\fp A_\fp$.  
Let $\Upsilon_\fp:\Spec(Z_\fp) \rarr \Spec(R_\fp)$ be the morphism induced by the canonical map $R_\fp \rarr Z_\fp$. Recall from Lemma \ref{ff_center_ext} that $R_\fp \subs Z_\fp$ is a finite extension so that $\Upsilon_\fp$ is surjective. We have
\begin{align*}
\Upsilon_{\fp }^{-1}(\fp _{\fp }) & = \lbrace \fQ \in \Spec(Z_{\fp }) \mid \fQ \cap R_{\fp } = \fp _{\fp } \rbrace \\&= \lbrace \fQ \in \Spec(Z_{\fp }) \mid \fp _{\fp } \subs \fQ \rbrace \\ &= \lbrace \fQ \in \Spec(Z_{\fp }) \mid \fp _{\fp }Z_{\fp } \subs \fQ \rbrace \\
& \simeq \Spec(Z(\fp)) \;.
\end{align*}
In the second equality we used the fact that $R_{\fp } \rarr Z_{\fp }$ is a finite morphism and $R_{\fp }$ is local with maximal ideal $\fp _{\fp }$. The identification with $\Spec(Z(\fp))$ is canonical since $Z(\fp) = Z_\fp/\fp_\fp Z_\fp$. The morphism $\Theta_\fp: \Spec(Z_{\fp }) \rarr \Spec(Z)$ induced by the localization map $Z \rarr Z_\fp$ is injective by \cite[Proposition 2.2(b)]{Eis-Commutative-Algebra}. We claim that this map induces $\Upsilon^{-1}_\fp(\fp_\fp) \simeq \Upsilon^{-1}(\fp)$. If $\fQ \in \Upsilon_{\fp }^{-1}(\fp _{\fp })$, then clearly $(\fQ \cap Z) \cap R = \fQ \cap R \subs R \cap \fp _{\fp } = \fp $ and therefore $\Theta_\fp$ induces an injective map $\Upsilon_{\fp }^{-1}(\fp _{\fp }) \rarr \Upsilon^{-1}(\fp )$. If $\fQ \in \Upsilon^{-1}(\fp )$, then, since $\fQ \cap R = \fp $, we have $\fQ \cap (R \setminus \fp ) = \emptyset$ so that $\fQ_{\fp } \in \Spec(Z_{\fp })$ and clearly $\fp _{\fp } \subs \fQ_{\fp }$, implying that $\fQ_{\fp } \in \Upsilon_{\fp }^{-1}(\fp _{\fp })$. The map $\Upsilon_{\fp }^{-1}(\fp _{\fp }) \rarr \Upsilon^{-1}(\fp )$ is thus bijective. Hence, we have a canonical bijection $\Spec(Z(\fp)) \simeq \Upsilon^{-1}(\fp)$. Now, recall from Lemma \ref{cliques_zero_dim} that $\Spec(Z(\fp)) \simeq \msf{Bl}(Z(\fp))$. 
\end{proof}

\subsection{Blocks via central characters} \hfill

\begin{tcolorbox}
We assume that $R$ is \textit{noetherian} and \textit{normal}, and that $A$ is a finite flat $R$-algebra with \textit{split} generic fiber $A^K$.
\end{tcolorbox}

Recall from Corollary \ref{normal_split_corollary} that the quotient map $A_\fp \twoheadrightarrow A(\fp)$ induces $\msf{Bl}(A_\fp) \simeq \msf{Bl}(A(\fp))$, so together with Theorem \ref{mueller_maintheorem} we have a canonical bijection
\begin{equation} \label{bl_ap_ups_bij}
\msf{Bl}(A_\fp) \simeq \Upsilon^{-1}(\fp) \;.
\end{equation}
Recall from \S\ref{blocks_of_localizations} that $\msf{Fam}_\fp(A^K)$ is the partition of $\Irr A^K$ induced by the blocks of $A_\fp$ and that we naturally have $\msf{Bl}(A_\fp) \simeq \msf{Fam}_\fp(A^K)$. Altogether, we now have canonical bijections
 \begin{equation} \label{pfam_fiber_bij}
 \msf{Fam}_\fp(A) \simeq \msf{Bl}(A_\fp) \simeq \Upsilon^{-1}(\fp) \simeq \msf{Bl}(A(\fp)) \;.
 \end{equation}
 Since $A$ has split generic fiber $A^K$, we have a central character $\Omega_S:\msf{Z}(A^K) \rarr K$ for every simple $A^K$-module $S$. Recall that $\Omega_S(z)$ is the scalar by which $z \in \msf{Z}(A^K)$ acts on $S$. Since $R$ is normal, the image of the restriction of $\Omega_S$ to $\msf{Z}(A) \subs \msf{Z}(A^K)$ is contained in $R \subs K$. We thus get a well-defined $R$-algebra morphism 
 \begin{equation}
\Omega_S':\msf{Z}(A) \rarr R \;.
\end{equation}
It is a classical fact that $S,T \in \Irr A^K$ lie in the same family if and only if $\Omega_S' = \Omega_T'$. We can thus label the central characters of $A^K$ as $\Omega_\mathcal{F}$ with $\mathcal{F}$ a family (block) of $A^K$. Using Theorem \ref{mueller_maintheorem} this description generalizes modulo $\fp$ so that we get an explicit description of the $\fp$-families, and thus of the block stratification. For $\fp \in \Spec(R)$ let
\begin{equation}
\Omega_S^\fp : \msf{Z}(A) \to R/\fp
\end{equation}
be the composition of $\Omega_S'$ with the quotient map $R \twoheadrightarrow R/\fp$.

\begin{thm} \label{maintheorem_blex_descr}
Under the bijection $\Upsilon^{-1}(\fp) \simeq \msf{Fam}_\fp(A)$ from (\ref{pfam_fiber_bij}) the $\fp$-family of a simple $A^K$-module $S$ corresponds to $\Ker \Omega_S^\fp$. Hence, two simple $A^K$-modules $S$ and $T$ lie in the same $\fp$-family if and only if $\Omega_S'(z) \equiv \Omega_T'(z) \ \msf{mod} \ \fp$ for all $z \in \msf{Z}(A)$. So, if $z_1,\ldots,z_n$ is an $R$-algebra generating system of $\msf{Z}(A)$ and $\mscr{F},\mscr{F}'$ are two distinct $A^K$-families, then the corresponding gluing locus is given by
\begin{equation}
\msf{Gl}_A(\leq \! \lbrace \mscr{F},\mscr{F}' \rbrace) = \msf{V}( \lbrace \Omega_\mscr{F}(z_i) - \Omega_{\mscr{F}'}(z_i) \mid i =1,\ldots,n \rbrace) \;.
\end{equation}
\end{thm}

\begin{proof}
Considering the explicit form of the bijection given in Theorem \ref{mueller_maintheorem} we see that the bijection (\ref{bl_ap_ups_bij}) maps a block idempotent $c$ of $A_\fp$ to the (by the theorem unique) maximal ideal $\fQ_c$ of $Z$ lying above $\fp$ and satisfying $c^\dagger \in (\fQ_c)_\fp$. Let $c_\fQ$ be the block idempotent of $A_\fp$ corresponding to $\fQ \in \Upsilon^{-1}(\fp)$.

For $S \in \Irr A^K$ let $\Omega_S^{\fp }:Z \rarr R/\fp $ be the composition of $\Omega_S'$ and the quotient morphism $R \rarr R/\fp $. It is clear that $\Ker(\Omega_S^{\fp })  \in \Upsilon^{-1}(\fp )$. Note that $\Omega_S'(z) \equiv \Omega_T'(z) \ \msf{mod} \ \fp$ for all $z \in \msf{Z}(A)$ if and only if $\Omega_S^{\fp } = \Omega_T^{\fp }$. We have an exact sequence
\[
0 \longrightarrow \Ker(\Omega_S') \longrightarrow Z \overset{\Omega_S'}{\longrightarrow} R \longrightarrow 0 
\]
of $R$-modules. Since $\Omega_S'$ is an $R$-algebra morphism, the canonical map $R \rarr Z$ is a section of $\Omega_S'$ and therefore $Z = R \oplus \Ker(\Omega_S')$ as $R$-modules. Similarly, we have $Z = R \oplus \Ker(\Omega_T')$. Since $\Ker(\Omega_S') \subs \Ker(\Omega_S^{\fp })$ and $\Ker(\Omega_T') \subs \Ker(\Omega_T^{\fp })$, this implies that $\Omega_S^\fp = \Omega_T^\fp$ if and only if $\Ker(\Omega_S^\fp) = \Ker(\Omega_T^\fp)$.  

Now, suppose that $\Ker(\Omega_S^{\fp }) = \Ker(\Omega_T^{\fp })$. Denote this common kernel by $\fQ$. Clearly, $\fQ \in \Upsilon^{-1}(\fp )$. We know that the corresponding block idempotent $c_{\fQ}$ of $A_{\fp }$ has the property that $c_{\fQ}^\dagger \in \fQ_{\fp }$. Since $\Ker(\Omega_S') \subs \Ker(\Omega_S^{\fp }) = \fQ = \Ker(\Omega_T^{\fp }) \sups \Ker(\Omega_T')$, this certainly implies that $c_{\fQ}^\dagger S= 0 = c_{\fQ}^\dagger T$. Hence, $S$ and $T$ lie in the same $\fp $-family.

Conversely, suppose that $S$ and $T$ lie in the same $\fp$-family. We can write the corresponding block idempotent of $A_\fp$ as $c_{\fQ}$ for some $\fQ \in \Upsilon^{-1}(\fp )$. By definition, $c_{\fQ}^\dagger S = 0 = c_{\fQ}^\dagger T$. We know that $c_{\fQ}^\dagger \in \fQ_{\fp }$ and $c_{\fQ} \notin \fQ_{\fp }$ and therefore $\Ker( (\Omega_S')_{\fp }) = \fQ_{\fp } = \Ker( (\Omega_T')_{\fp }) $. Hence, $\fQ \subs \Ker(\Omega_S') \subs \Ker(\Omega_S^{\fp })$ and $\fQ \subs \Ker(\Omega_T') \subs \Ker(\Omega_T^{\fp })$. Since $\fQ, \Ker(\Omega_S^{\fp }), \Ker(\Omega_T^{\fp }) \in \Upsilon^{-1}(\fp )$ and all prime ideals in $\Upsilon^{-1}(\fp )$ are incomparable, we thus conclude that $\Ker(\Omega_S^{\fp }) = \Ker(\Omega_T^{\fp })$. 

The equation for the gluing locus is now clear.
\end{proof}

\section{Blocks and decomposition matrices} \label{blocks_and_dec_maps}

To obtain information about the actual members of the $A(\fp)$-families we use decomposition maps as introduced by Geck and Rouquier \cite{Geck-Rouquier-Dec} (see also \cite{Geck-Pfeiffer} and \cite{Thiel-Dec}). For a ring $A$ we denote by $\msf{G}_0(A) \dopgleich  \msf{K}_0(A\tn{-}\msf{mod})$ the \word{Grothendieck group} and by $\msf{K}_0(A) \dopgleich \msf{K}_0(A\tn{-}\msf{proj})$ the \word{projective class group}. In case $A$ is semiperfect (e.g., artinian), $\msf{K}_0(A)$ is the free abelian group with basis the isomorphism classes of the projective indecomposable modules. In case $A$ is artinian, $\msf{G}_0(A)$ is the free abelian group with basis the isomorphism classes of simple modules and $\msf{K}_0(A) \simeq \msf{G}_0(A)$ mapping $P$ to $\Hd(P)$.

For the theory of decomposition maps we need the following (standard) assumption:

\begin{tcolorbox}
$A$ is finite free with split generic fiber and for any non-zero  $\fp \in \Spec(R)$ there is a discrete valuation ring $\sO$ with maximal ideal $\fm$  in $K$ dominating $R_\fp$ such that the canonical map $\msf{G}_0(A(\fp)) \rarr \msf{G}_0(A^\sO(\fm))$ of Grothendieck groups is an isomorphism.
\end{tcolorbox}

We call  a ring $\sO$ as above a \word{perfect $A$-gate} in $\fp$. We refer to \cite{Thiel-Dec} for more details. The following lemma lists two standard situations in which the above assumptions hold. Part \ref{dec_theory_assumptions_lemma:dedekind} is obvious and part \ref{dec_theory_assumptions_lemma:noeth} was proven in \cite[Theorem 1.22]{Thiel-Dec}.

\begin{lem} \label{dec_theory_assumptions_lemma}
A finite free $R$-algebra $A$ with split generic fiber satisfies the above assumptions in the following two cases:
\begin{enum_thm}
\item \label{dec_theory_assumptions_lemma:dedekind} $R$ is a Dedekind domain.
\item \label{dec_theory_assumptions_lemma:noeth} $R$ is noetherian and $A$ has split fibers.
\end{enum_thm}
\end{lem}

If $\sO$ is a perfect $A$-gate in $\fp$, then there is a group morphism 
\begin{equation}
\msf{d}_A^{\fp,\sO}:\msf{G}_0(A^K) \rarr \msf{G}_0(A(\fp))
\end{equation}
between Grothendieck groups generalizing reduction modulo $\fp$. In case $R$ is normal, it was proven by Geck and Rouquier \cite{Geck-Rouquier-Dec} that this map is independent of the choice of $\sO$ and in this case we just write $\msf{d}_A^\fp$. We note that in case $R$ is noetherian and $A$ has split fibers, any decomposition map in the sense of Geck and Rouquier can be realized by a perfect $A$-gate, see \cite[Theorem 1.22]{Thiel-Dec}.

\subsection{Brauer reciprocity}

An important tool for relating decomposition maps and blocks is \word{Brauer reciprocity} which we prove in Theorem \ref{brauer_rec} below in our general setup (this was known to hold before only in special settings). Recall that the \word{intertwining form} for a finite-dimensional algebra $B$ over a field $F$ is the $\bbZ$-linear pairing $\langle \cdot, \cdot \rangle_B:\msf{K}_0(B) \times \msf{G}_0(B) \rarr \bbZ$ uniquely defined by 
\begin{equation} \label{intertwining_form_def}
\langle \lbrack P \rbrack, \lbrack V \rbrack \rangle \dopgleich \dim_F \Hom_B(P,V)
\end{equation}
for a finite-dimensional projective $B$-module $P$ and a finite-dimensional $B$-module $V$, see \cite[\S2]{Geck-Rouquier-Dec}. Here, $\msf{K}_0(B)$ is the zeroth $\msf{K}$-group of the category of finite-dimensional projective $B$-modules. The intertwining form is always non-degenerate, see Lemma \ref{intertwining_nondeg}. Due to the non-degeneracy of $\langle \cdot,\cdot \rangle_{A^K}$ there is at most one adjoint 
\begin{equation}
\msf{e}_A^{\fp,\sO}:\msf{K}_0(A(\fp)) \rarr \msf{K}_0(A^K)
\end{equation}
 of $\msf{d}_A^{\fp,\sO}:\msf{G}_0(A^K) \rarr \msf{G}_0(A(\fp))$ with respect to $\langle \cdot,\cdot \rangle_{A(\fp)}$, characterized by the relation
\begin{equation} \label{brauer_rec_adjoint_rel}
\langle \msf{e}_A^{\fp,\sO}(\lbrack \ol{P} \rbrack), \lbrack V \rbrack \rangle_{A^K} = \langle \lbrack \ol{P} \rbrack, \msf{d}_A^{\fp,\sO}(\lbrack V \rbrack) \rangle_{A(\fp)} \;.
\end{equation}
for all finitely generated $A^K$-modules $V$ and all finitely generated projective $A(\fp)$-modules $\ol{P}$, see Lemma \ref{intertwining_nondeg}. Brauer reciprocity is about the existence of this adjoint.

\begin{thm} \label{brauer_rec}
The (unique) adjoint $\msf{e}_A^{\fp,\sO}$ of $\msf{d}_A^{\fp,\sO}$ exists. Moreover, the diagram
\begin{equation} \label{brauer_rec_diagram}
\begin{tikzcd}
\msf{K}_0(A^K) \ar{r}{\msf{c}_{A^K}} & \msf{G}_0(A^K) \ar{d}{\msf{d}_A^{\fp,\sO}} \\
\msf{K}_0(A(\fp)) \ar{r}[swap]{\msf{c}_{A(\fp)}} \ar{u}{\msf{e}_A^{\fp,\sO}} & \msf{G}_0(A(\fp)) 
\end{tikzcd}
\end{equation}
commutes, where the horizontal morphisms are the canonical ones (\word{Cartan maps}) mapping a class $\lbrack P \rbrack$ of a projective module $P$ to its class $\lbrack P \rbrack$ in the Grothendieck group. If $R$ is normal, the morphism $\msf{e}_A^{\fp,\sO}$ does not depend on the choice of $\sO$ and we denote it by $\msf{e}_A^\fp$.
\end{thm}

\begin{proof}
Since $\langle \cdot, \cdot \rangle_{A^K}$ is non-degenerate by Lemma \ref{intertwining_nondeg}, it follows that $\msf{d}_A^{\fp,\sO}$ has at most one adjoint $\msf{e}_A^{\fp,\sO}$, characterized by equation (\ref{brauer_rec_adjoint_rel}), see \cite[Satz 78.1]{Scheja-Storch-Alg-2}. By assumption  there is a perfect $A$-gate $\sO$ in $\fp$. Let $\fm$ be the maximal ideal of $\sO$. Since $A^K$ splits by assumption, Corollary \ref{algebra_dvr_semiperfect} implies that $A^\sO$ is semiperfect. The morphism $\msf{K}_0(A^\sO) \rarr \msf{K}_0(A^\sO(\fm))$ induced by the quotient map $A^\sO \twoheadrightarrow A^\sO(\fm)$ is thus an isomorphism by lifting of idempotents. Furthermore, by assumption the morphism $\msf{d}_A^{\fp,\fm}:\msf{G}_0(A(\fp)) \rarr \msf{G}_0(A^\sO(\fm))$ is an isomorphism and then the proof of Theorem \ref{faithfully_flat_ext_block_bij} shows that the canonical morphism $\msf{e}_A^{\fp,\fm}:\msf{K}_0(A(\fp)) \rarr \msf{K}_0(A^\sO(\fm))$ is also an isomorphism. We can thus define a morphism $\msf{e}_A^{\fp,\sO}:\msf{K}_0(A(\fp)) \rarr \msf{K}_0(A^K)$ as the following composition
\begin{equation}
\begin{tikzcd}
\msf{K}_0(A(\fp)) \arrow{r}{\simeq} \arrow[bend right=15]{rrr}[swap]{\msf{e}_A^{\fp,\sO}} & \msf{K}_0(A^\sO(\fm)) \arrow{r}{\simeq} & \msf{K}_0(A^\sO) \arrow{r} & \msf{K}_0(A^K)
\end{tikzcd}
\end{equation}

We will now show that $\msf{e}_A^{\fp,\sO}$ is indeed an adjoint of $\msf{d}_A^{\fp,\sO}$. The arguments in the proof of \cite[18.9]{CR-Methods-1} can, with some refinements, be transferred to our more general situation and this is what we will do. Let $\ol{P}$ be a finitely generated projective $A(\fp)$-module and let $V$ be a finitely generated $A^K$-module. Since $\msf{K}_0(A^\sO) \simeq \msf{K}_0(A^\sO(\fm))$,  there exists a finitely generated projective $A^{\sO}$-module $P$ such that $(\msf{e}_{A}^{\fp,\fm})^{-1}(\lbrack P/\fm P \rbrack) = \lbrack \ol{P} \rbrack$ and then we have $\msf{e}_A^{\fp,\sO}(\lbrack \ol{P} \rbrack) = \lbrack P^K \rbrack$. Let $\wt{V}$ be an $A^{\sO}$-lattice in $V$. Then by definition of $\msf{d}_A^{\fp,\sO}$, see \cite[Corollary 1.14]{Thiel-Dec}, we have $\msf{d}_A^{\fp,\sO}(\lbrack V \rbrack) = (\msf{d}_A^{\fp,\fm})^{-1}(\lbrack \wt{V}(\fm) \rbrack)$. We denote by $\ol{V}$ a representative of $\msf{d}_A^{\fp,\sO}(\lbrack V \rbrack)$. Since $P$ is a finitely generated projective $A^{\sO}$-module, we can write $P \oplus Q = (A^{\sO})^n$ for some finitely generated projective $A^\sO$-module $Q$ and some $n \in \bbN$. Since $\Hom_{A^{\sO}}$ is additive, we get 
\begin{align*}
\Hom_{A^{\sO}}(P, \wt{V}) \oplus \Hom_{A^{\sO}}(Q,\wt{V}) & = \Hom_{A^{\sO}}(P \oplus Q, \wt{V}) = \Hom_{A^{\sO}}( (A^{\sO})^n, \wt{V}) \\ &= ( \Hom_{A^{\sO}}( A^{\sO}, \wt{V}) )^n \simeq \wt{V}^n \;.
\end{align*}
This shows that $\Hom_{A^{\sO}}(P,\wt{V})$ is a direct summand of $\wt{V}^n$ and as $\wt{V}^n$ is $\sO$-free, we conclude that $\Hom_{A^{\sO}}(P,\wt{V})$ is $\sO$-projective and thus even $\sO$-free since $\sO$ is a discrete valuation ring. Since $P$ is a finitely generated projective $A^\sO$-module, it follows from Lemma \ref{base_ring_change_of_hom} that there is a canonical $K$-vector space isomorphism 
\[
K \otimes_{\sO} \Hom_{A^{\sO}}(P,\wt{V}) \simeq \Hom_{A^K}(P^K,V)
\]
and a canonical $\msf{k}(\fm)$-vector space isomorphism
\[
\msf{k}(\fm) \otimes_{\sO} \Hom_{A^{\sO}}(P,\wt{V}) \simeq \Hom_{A^{\sO}(\fm)}(P/ \fm P, \wt{V}/\fm \wt{V}) \;.
\]
Combining all results and the fact that both $\msf{e}_A^{\fp,\sO}$ and $\msf{d}_A^{\fp,\sO}$ preserve dimensions by construction, we can now conclude that
\begin{align*}
\langle \msf{e}_A^{\fp,\sO}(\lbrack \ol{P} \rbrack), \lbrack V \rbrack \rangle_{A^K} &= \dim_K \Hom_{A^K}( P^K, V) = \dim_{\sO} \Hom_{A^{\sO}}( P, \wt{V}) \\ &= \dim_{\msf{k}(\fm)} \Hom_{A^{\sO}(\fm)}( P/\fm P, \wt{V}/\fm \wt{V}) = \dim_{\msf{k}(\fp)} \Hom_{A(\fp)}( \ol{P}, \ol{V}) \\ &= \langle \lbrack \ol{P} \rbrack, \msf{d}_A^{\fp,\sO}(\lbrack V \rbrack) \rangle_{A(\fp)} \;.
\end{align*}

Proving the commutativity of diagram (\ref{brauer_rec_diagram}) amounts to proving that $\msf{c}_{A(\fp)}( \lbrack \ol{P} \rbrack ) = \msf{d}_A^{\fp,\sO} \circ \msf{c}_{A^K} \circ \msf{e}_A^{\fp,\sO} ( \lbrack \ol{P} \rbrack)$ for every finitely generated projective $A(\fp)$-module $\ol{P}$. To prove this, note that the diagram
\[
\begin{tikzcd}
\msf{K}_0(A^{\sO}(\fm)) \arrow{r}{\msf{c}_{A^{\sO}(\fm)}} & \msf{G}_0(A^{\sO}(\fm)) \\
\msf{K}_0(A(\fp)) \arrow{r}[swap]{\msf{c}_{A(\fp)}} \arrow{u}{\msf{e}_A^{\fp,\fm}} & \msf{G}_0(A(\fp)) \arrow{u}[swap]{\msf{d}_A^{\fp,\fm}}
\end{tikzcd}
\]
commutes. As above we know that there exists a finitely generated projective $A^{\sO}$-module $P$ such that $(\msf{e}_A^{\fp,\fm})^{-1}(\lbrack P/\fm P \rbrack) = \lbrack \ol{P} \rbrack$ and $\msf{e}_A^{\fp,\sO}(\lbrack \ol{P} \rbrack) = \lbrack P^K \rbrack$. Since $P$ is a finitely generated projective $A^{\sO}$-module and $A$ is a finite $\sO$-module, it follows that $P$ is also a finitely generated projective $\sO$-module. As $\sO$ is a discrete valuation ring, we conclude that $P$ is actually $\sO$-free of finite rank. Hence, $P$ is an $A^{\sO}$-lattice in $P^K$ and therefore 
\begin{align*}
\msf{d}_A^{\fp,\sO} \circ \msf{c}_{A^K} \circ \msf{e}_A^{\fp,\sO} ( \lbrack \ol{P} \rbrack) &=  \msf{d}_A^{\fp,\sO} ( \lbrack P^K \rbrack) = (\msf{d}_{A}^{\fp,\fm})^{-1} ( \lbrack P/\fm P \rbrack) \\& = (\msf{d}_{A}^{\fp,\fm})^{-1} \circ \msf{c}_{A(\fm)}( \lbrack P/\fm P \rbrack)  \\& = \msf{c}_{A(\fp)} \circ (\msf{e}_{A}^{\fp,\fm})^{-1}( \lbrack P/\fm P \rbrack) = \msf{c}_{A(\fp)} ( \lbrack \ol{P} \rbrack). 
\end{align*}
If $R$ is normal, then the independence of $\msf{e}_A^{\fp,\sO}$ from the choice of $\sO$ follows from the independence of $\msf{d}_A^{\fp,\sO}$ from the choice of $\sO$ and the fact that $\msf{d}_A^{\fp,\sO}$ has at most one adjoint.
\end{proof}

\subsection{Preservation of simple modules vs. preservation of blocks} \label{preservation}

In \cite{Thiel-Dec} we studied the set
\begin{equation}
\msf{DecGen}(A) \dopgleich \left\lbrace \fp \in \Spec(R) \mid  \msf{d}_A^{\fp,\sO} \tn{ is trivial for any }A\tn{-gate in }\fp \right\rbrace \;.
\end{equation}
where $\msf{d}_A^{\fp,\sO}$ being \word{trivial} means that it induces a bijection between simple modules. We have proven in \cite[Theorem 2.3]{Thiel-Dec} that $\msf{DecGen}(A)$ is open if $R$ is noetherian and $A$ has split fibers. Brauer reciprocity implies that $\msf{e}_A^{\fp,\sO}$ is trivial if and only if $\msf{d}_A^{\fp,\sO}$ is trivial, so we deduce that the locus of all $\fp$ such that $\msf{e}_A^{\fp,\sO}$ is trivial for any $\sO$ is an open subset of $\msf{Spec}(R)$. 

If $\fp \in \msf{DecGen}(A)$, then the simple modules of $A^K$ and $A(\fp)$ are ``essentially the same'', in particular their dimensions are the same. This is why explicit knowledge about $\msf{DecGen}(A)$ is quite helpful to understand the representation theory of the fibers of $A$, see \cite{Thiel-Dec}. So far, we do not have an explicit description of $\msf{DecGen}(A)$, however. Brauer reciprocity enables us to prove the following relation between decomposition maps and blocks.

\begin{thm} \label{decgen_blgen_inclusion_main_theorem}
We have an inclusion
\begin{equation} \label{decgen_blgen_inclusion}
\msf{DecGen}(A) \subs \msf{BlGen}(A) \;.
\end{equation}
\end{thm}

\begin{proof}
Let $\fp \in \Spec(R)$ be non-zero. By assumption there is a perfect $A$-gate $\sO$ in $\fp$. If $\fp \in \msf{DecGen}(A)$, then by definition $\msf{d}_A^{\fp,\sO}$ is trivial, so the matrix $\msf{D}_A^{\fp,\sO}$ of this morphism in bases given by isomorphism classes of simple modules of $A^K$ and $A(\fp)$, respectively, is equal to the identity matrix when ordering the bases appropriately. It now follows from Brauer reciprocity, Theorem \ref{brauer_rec}, that $\msf{C}_{A(\fp)} = \msf{C}_{A^K}$ in appropriate bases, where $\msf{C}_{A(\fp)}$ is the matrix of the Cartan map $\msf{c}_{A(\fp)}$ and $\msf{C}_{A^K}$ is the matrix of the Cartan map $\msf{c}_{A^K}$. Due to the linkage relation explained in \S\ref{notations}, the families of $A^K$ and of $A(\fp)$ are determined by the respective Cartan matrices. Since $\msf{C}_{A(\fp)} = \msf{C}_{A^K}$, it follows that $\#\msf{Bl}(A(\fp)) = \#\msf{Bl}(A^K)$, so $\fp \in \msf{BlGen}(A)$.
\end{proof}

\begin{rem}
Suppose that $A$  has split fibers and that $R$ is noetherian. Then the fact that $\#\msf{Bl}(A(\fp)) = \#\msf{Bl}(Z(\fp))$ by Theorem \ref{mueller_maintheorem} together with the Lemma \ref{num_blocks_rad_dim} yields the following equivalence:
\begin{equation} \label{blgen_equivalent_condition_with_rad}
\fp \in \msf{BlGen}(A) \Longleftrightarrow \dim_K (Z^K + \Rad(A^K)) = \dim_{\msf{k}(\fp)} (Z(\fp) + \Rad(A(\fp))) \;.
\end{equation}
Let $\sO$ be a perfect $A$-gate in $\fp$. This exists by Lemma \ref{dec_theory_assumptions_lemma}\ref{dec_theory_assumptions_lemma:noeth}.  Suppose that $\fp \in \msf{DecGen}(A)$. In \cite[Theorem 2.2]{Thiel-Dec} we have proven that this implies that \[ \dim_K \Rad(A^K) = \dim_{\msf{k}(\fp)} \Rad(A(\fp))\;. \] Let $X \dopgleich Z + J$, where $J \dopgleich \Rad(A^K) \cap A^\sO$. The arguments in \cite{Thiel-Dec} show that $X$ is an $A^\sO$-lattice of $Z^K + \Rad(A^K)$ and that the reduction in the maximal ideal $\fm$ of $\sO$ is equal to $Z^\sO(\fm) + \Rad(A^\sO(\fm))$. We thus have $\dim_K(Z^K + \Rad(A^K)) = \dim_{\msf{k}(\fm)}(Z^\sO(\fm) + \Rad(A^\sO(\fm)))$. Since $A(\fp)$ splits, the $\msf{k}(\fm)$-dimension of $Z^\sO(\fm) + \Rad(A^\sO(\fm))$ is equal to the $\msf{k}(\fp)$-dimension of $Z(\fp) + \Rad(A(\fp))$. Hence, we have $\fp \in \msf{BlGen}(A)$ by (\ref{blgen_equivalent_condition_with_rad}). This yields another proof of the inclusion $\msf{DecGen}(A) \subs \msf{BlGen}(A)$ in case $A$ has split fibers.
\end{rem}

\begin{ex} \label{decgen_eq_blgen_counterex}
The following example due to C. Bonnafé shows that in the generality of Theorem \ref{decgen_blgen_inclusion_main_theorem} we do not have equality in (\ref{decgen_blgen_inclusion}). Let $R$ be a discrete valuation ring with fraction field $K$ and uniformizer $\pi$, i.e., $\fp \dopgleich (\pi)$ is the maximal ideal of $R$. Denote by $k \dopgleich R/\fp$ the residue field in $\fp$. Let
\[
A \dopgleich \left\lbrace \begin{pmatrix} a & b \\ c & d \end{pmatrix} \in \Mat_2(R) \mid b,c \in \fp \right\rbrace \;.
\]
This is an $R$-subalgebra of $\Mat_2(R)$ and it is $R$-free with basis
\begin{equation} \label{decgen_eq_blgen_counterex_gens}
e \dopgleich E_{11} , \; f \dopgleich E_{22} , \; x \dopgleich \pi E_{12} , \; y \dopgleich \pi E_{21},
\end{equation}
where $E_{ij} = (\delta_{i,k}\delta_{j,l})_{kl}$ is the elementary matrix. Clearly, $A^K = \Mat_{2}(K)$, so the generic fiber of $A$ is split semisimple. In particular, $A^K$ has just one block, and this block contains just one simple module we denote by $S$. Now, consider the specialization $\ol{A} \dopgleich A(\fp) = A/\fp A$. We know from Corollary \ref{algebra_dvr_semiperfect} that the quotient map $A \twoheadrightarrow \ol{A}$, $a \mapsto \ol{a}$, is block bijective, so we must have $\#\msf{Bl}(A(\fp)) \leq \msf{Bl}(A^K)$ and therefore $\# \msf{Bl}(A(\fp)) = 1$, so $\fp \in \msf{BlGen}(A)$. Let $\ol{J}$ be the $k$-subspace of $\ol{A}$ generated by $\ol{x}$ and $\ol{y}$. This is in fact a two-sided ideal of $\ol{A}$ since it is stable under multiplication by the generators (\ref{decgen_eq_blgen_counterex_gens}). Moreover, we have $\ol{x}^2 = 0 = \ol{y}^2$, so $\ol{J}$ is a nilpotent ideal of $\ol{A}$. Hence, $\dim_k \Rad(\ol{A}) \geq 2$. The number of simple modules of $\ol{A}$ is by \cite[Theorem 7.17]{Lam-First-Course-91} equal to $\dim_k \ol{A}/(\Rad(\ol{A}) + \lbrack \ol{A},\ol{A} \rbrack)$, so $\#\Irr \ol{A} \leq 2$ since $\dim_k \ol{A} = \dim_K A^K = 4$. The two elements $\ol{e}$ and $\ol{f}$ are orthogonal idempotents and so the constituents of the two $\ol{A}$-modules $\ol{A} \ol{e}$ and $\ol{A}\ol{f}$ are non-isomorphic. So, we have $\# \Irr \ol{A} \geq 2$ and due to the aforementioned we conclude that $\#\Irr \ol{A} = 2$. Let $\ol{S}_1$ and $\ol{S}_2$ be these two simple modules. Since $R$ is a discrete valuation ring, reduction modulo $\fp$ yields the well-defined decomposition map $\msf{d}_A^\fp:\msf{G}_0(A^K) \rarr \msf{G}_0(A(\fp))$, see \cite[Corollary 1.14]{Thiel-Dec}. It is an elementary fact that the all simple $\ol{A}$-modules must be constituents of $\msf{d}_A^\fp(\lbrack S \rbrack) = \lbrack S/\fp S \rbrack$. Since $\dim_K S = 2$, the only possibility is that $\msf{d}_A^\fp(\lbrack S \rbrack) = \lbrack \ol{S}_1 \rbrack + \lbrack \ol{S}_2 \rbrack$ and $\dim_k \ol{S}_i = 1$. In particular, $\fp \notin \msf{DecGen}(A)$, so $\fp \in \msf{BlGen}(A) \setminus \msf{DecGen}(A)$. Finally, we note that $\ol{A}$ also splits since $\#\Irr(\ol{A}) = 2$ implies by the above formula that $\dim_k \Rad(\ol{A}) = 2$ and we have $\dim_k \ol{A} = \dim_k \Rad(\ol{A}) + \sum_{i=1}^2 (\dim_k \ol{S}_i)^2$, so $\ol{A}$ is split by \cite[Corollary 7.8]{Lam-First-Course-91}.  
\end{ex}

\begin{lem} \label{blgen_decgen_complement_diagonal}
Assume that the $A^K$-families are singletons, and that $\#\Irr A(\fp) \leq \# \Irr A^K$ for all $\fp \in \Spec(R)$. Then
\[
\msf{BlGen}(A) \setminus \msf{DecGen}(A) = \lbrace \fp \in \Spec(R) \mid \msf{D}_A^\fp \tn{ is diagonal but not the identity } \rbrace \;.
\]
\end{lem}

\begin{proof}
Since the $A^K$-families are singletons, we have $\#\msf{Irr}(A^K) = \#\msf{Bl}(A^K)$. We clearly have $\#\Irr A(\fp) \geq \#\msf{Bl}(A(\fp))$ for all $\fp \in \Spec(R)$. Assume that $\fp \in \msf{BlGen}(A)$. Then we have $\#\Irr A(\fp) \geq \#\msf{Irr}(A^K)$, so $\#\Irr A(\fp)) = \#\Irr A^K$ by our assumption. Hence, the decomposition matrix $\msf{D}_A^\fp$ is quadratic. By Theorem \ref{dec_block_compat} the $\fp$-families are equal to the Brauer $\fp$-families. Since $\fp \in \msf{BlGen}(A)$ and the $A^K$-families are singletons, it follows that $\msf{D}_A^\fp$ is a diagonal matrix. The claim is now obvious. 
\end{proof}

\begin{lem} \label{cellular_alg_blgen_decgen}
Let $R$ be a noetherian integral domain with fraction field $K$ and let $A$ be a cellular $R$-algebra of finite dimension such that $A^K$ is semisimple. Then $\msf{DecGen}(A) = \msf{BlGen}(A)$. 
\end{lem}

\begin{proof}
First of all, specializations of $A$ are again cellular by \cite[1.8]{Graham-Lehrer-Cellular}. Moreover, it follows from \cite[Proposition 3.2]{Graham-Lehrer-Cellular} that $A$ has split fibers, so $A$ satisfies Lemma \ref{dec_theory_assumptions_lemma}\ref{dec_theory_assumptions_lemma:noeth} and therefore our basic assumption in this paragraph. Let $\Lambda$ be the poset of the cellular structure of $A^K$. Since $A^K$ is semisimple, each cell module $M_\lambda$ has simple head $S_\lambda$ and $\#\Irr A^K = \#\Lambda$. 
Let $\fp \in \Spec(R)$. The poset for the cellular structure of $A(\fp)$ is again $\Lambda$. Denote by $M_\lambda^\fp$ the corresponding cell modules of $A(\fp)$. There is a subset $\Lambda'$ of $\Lambda$ such that $M_\lambda^\fp$ has simple head $S_\lambda^\fp$ for all $\lambda \in \Lambda'$ and that these heads are precisely the simple $A(\fp)$-modules. In particular, we have $\#\Irr A(\fp) \leq \#\Irr A^K$. Now, assume that $\fp \in \msf{BlGen}(A)$. By Lemma \ref{blgen_decgen_complement_diagonal} we just need to show that the decomposition matrix $\msf{D}_A^\fp$, which is square by the proof of Lemma \ref{blgen_decgen_complement_diagonal}, cannot be a non-identity diagonal matrix. By \cite[Proposition 3.6]{Graham-Lehrer-Cellular} we know that $\lbrack M_\lambda: S_\lambda \rbrack = 1$ and $\lbrack M_\lambda^\fp:S_\lambda^\fp \rbrack = 1$. By construction, it is clear that $\msf{d}_A^\fp(\lbrack M_\lambda \rbrack) = \lbrack M_\lambda^\fp \rbrack$. Hence, if $\msf{d}_A^\fp(\lbrack S_\lambda \rbrack) = n_\lambda \lbrack S_\lambda^\fp \rbrack$, we have $n_\lambda = \lbrack M_\lambda^\fp:S_\lambda^\fp \rbrack = 1$. Hence, $\msf{D}_A^\fp$ is the identity matrix, so $\fp \in \msf{BlGen}(A)$.
\end{proof}

\subsection{The Brauer graph} \label{brauer_graph}

Geck and Pfeiffer \cite{Geck-Pfeiffer} have introduced the so-called Brauer $\fp$-graph of $A$ in our general context but assuming that $A^K$ is semisimple so that the $A^K$-families are singletons. For general $A$ this definition seems not to be the correct one. We introduce the following generalization of this concept.

\begin{defn}
Suppose that $R$ is normal so that we have unique decomposition maps. The \word{Brauer $\fp$-graph} of $A$ is the graph with vertices the simple $A^K$-modules and an edge between $S$ and $T$ if and only if in the $A^K$-family of $S$ there is some $S'$ and in the $A^K$-family of $T$ there is some $T'$ such that $\msf{d}_A^\fp(\lbrack S' \rbrack)$ and $\msf{d}_A^\fp(\lbrack T' \rbrack)$ have a common constituent. The connected components of this graph are called the \word{Brauer $\fp$-families} of $A$. 
\end{defn}

If the $A^K$-families are singletons, we have an edge between $S$ and $T$ if and only if $\msf{d}_A^\fp(\lbrack S \rbrack)$ and $\msf{d}_A^\fp(\lbrack T \rbrack)$ have a common constituent, so this indeed generalizes the Brauer $\fp$-graph from \cite{Geck-Pfeiffer} for $A^K$ semisimple. Our final theorem shows that decomposition maps are compatible with $\fp$-families and $A(\fp)$-families, and relates the Brauer $\fp$-families to the $\fp$-families.

\begin{thm} \label{dec_block_compat}
Assume that $R$ is normal. The following holds:
\begin{enum_thm}
\item \label{dec_block_compat:a} A finite-dimensional $A^K$-module $V$ belongs to a $\fp$-block of $A$ if and only if $\msf{d}_A^\fp(\lbrack V \rbrack)$ belongs to a block of $A(\fp)$.
\item \label{dec_block_compat:b} Two finite-dimensional $A^K$-modules $V$ and $W$ lie in the same $\fp$-block if and only if $\msf{d}_A^\fp(\lbrack V \rbrack)$ and $\msf{d}_A^\fp(\lbrack W \rbrack)$ lie in the same block of $A(\fp)$.
\item If $\mathcal{F} \in \msf{Fam}_\fp(A)$ is a $\fp$-family, then 
\[
\msf{d}_A^\fp(\mathcal{F}) \dopgleich \lbrace T \mid T \tn{ is a constituent of } \msf{d}_A^\fp(\lbrack S \rbrack) \tn{ for some } S \in \mathcal{F} \rbrace
\]
is a family of $A(\fp)$, and all families of $A(\fp)$ are obtained in this way.
\item The Brauer $\fp$-families are equal to the $\fp$-families.
\end{enum_thm}
\end{thm}

\begin{proof}
\hfill

\begin{enum_proof}
\item By assumption there is a perfect $A$-gate $\sO$ in $\fp$. Let $\fm$ be the maximal ideal of $\sO$. We have the following commutative diagram of canonical morphisms which are all idempotent stable:
\begin{equation}
\begin{tikzcd}
A_\fp \arrow[hookrightarrow]{r} \arrow[twoheadrightarrow]{d} & A^\sO \arrow[hookrightarrow]{r} \arrow[twoheadrightarrow]{d} & A^K \\
A(\fp) \arrow[hookrightarrow]{r} & A^\sO(\fm)
\end{tikzcd}
\end{equation}
Since $R$ is assumed to be normal, it follows from Proposition \ref{unibranched_block_bijective} that  $A_\fp \twoheadrightarrow A(\fp)$ is block bijective. By assumption  the morphism $\msf{d}_A^{\fp,\fm}: \rG_0(A(\fp)) \rarr \rG_0(A^\sO(\fm))$  is an isomorphism and therefore $A(\fp) \hookrightarrow A^\sO(\fm)$ is block bijective by Theorem \ref{faithfully_flat_ext_block_bij}. Furthermore, by assumption the generic fiber $A^K$ is split and therefore $A^\sO \twoheadrightarrow A^\sO(\fm)$ is block bijective by Corollary \ref{algebra_dvr_semiperfect}. Because of (\ref{get_more_blocks_equation}) it thus follows that $A_\fp \hookrightarrow A^\sO$ is block bijective.

Now, let $V$ be a finite-dimensional $A^K$-module and let $\wt{V}$ be an $A^\sO$-lattice of $V$. Suppose that $V$ belongs to an $A_\fp$-block of $A^K$. Since $A_\fp \hookrightarrow A^\sO$ is block bijective, the $A_\fp$-blocks of $A^K$ coincide with the $A^\sO$-blocks of $A^K$ and therefore $V$ belongs to an $A^\sO$-block of $A^K$. Since $\wt{V}$ is $\sO$-free, it follows from Lemma \ref{block_extension_compatibility} that $\wt{V}$ belongs to a block of $A^\sO$. Again by Lemma \ref{block_extension_compatibility} and the fact that $A^\sO \twoheadrightarrow A^\sO(\fm)$ is block bijective, it follows that $\wt{V}/\fm \wt{V}$ belongs to a block of $A^\sO(\fm)$. Since $A(\fp) \hookrightarrow A^\sO(\fm)$ is block bijective, Lemma \ref{block_extension_compatibility} shows that $\msf{d}_A^\fp(\lbrack V \rbrack)$ belongs to a block of $A(\fp)$. 

Conversely, suppose that $\msf{d}_A^\fp(\lbrack V \rbrack)$ belongs to a block of $A(\fp)$. Then $\wt{V}/\fm \wt{V}$ belongs to a block of $A^\sO(\fm)$ and therefore $\wt{V}$ belongs to a block of $A^\sO$ by Lemma \ref{block_extension_compatibility}. But then $V$ belongs to an $A^\sO$-block of $A^K$ and thus to an $A_\fp$-block of $A^K$ by Lemma \ref{block_extension_compatibility}.

\item This follows now from part \ref{dec_block_compat:a}.

\item Fix a $\fp$-family $\mathcal{F}$ of $A^K$. If $S \in \mathcal{F}$, then $\msf{d}_A^\fp(\lbrack S \rbrack)$ belongs to an $A(\fp)$-block by \ref{dec_block_compat:a} and therefore all constituents of $\msf{d}_A^\fp(\lbrack S \rbrack)$ belong to a fixed family $\ol{\mathcal{F}}_S$. If $S' \in \mathcal{F}$ is another simple module, then by \ref{dec_block_compat:b} the constituents of $\msf{d}_A^\fp(\lbrack S' \rbrack)$ also lie in $\ol{\mathcal{F}}_S$. Hence, $\msf{d}_A^\fp(\mathcal{F})$ is contained in a fixed $A(\fp)$-family $\ol{\mathcal{F}}$. Let $T \in \ol{\mathcal{F}}$ be arbitrary. Due to the properties of decomposition maps there is some $S \in \Irr A^K$ such that $T $ is a constituent of $\msf{d}_A^\fp(\lbrack S \rbrack)$. Since $T$ and $\msf{d}_A^\fp(\lbrack S \rbrack)$ lie in the same $A(\fp)$-block by \ref{dec_block_compat:a} and \ref{dec_block_compat:b}, we must have $S \in \mathcal{F}$ by \ref{dec_block_compat:b}. Hence, $\ol{\mathcal{F}} = \msf{d}_A^\fp(\mathcal{F})$ is an $A(\fp)$-family. Since every simple $A(\fp)$-module is a constituent of $\msf{d}_A^\fp(\lbrack S \rbrack)$ for some simple $A^K$-module $S$, it is clear that any $A(\fp)$-family is of the form $\msf{d}_A^\fp(\mathcal{F})$ for a $\fp$-family $\mathcal{F}$.

\item Let $S$ and $T$ be simple $A^K$-modules contained in the same Brauer $\fp$-family, i.e., in the $A^K$-family of $S$ there is some $S'$ and in the $A^K$-family of $T$ there is some $T'$ such that $\msf{d}_A^\fp(\lbrack S' \rbrack)$ and $\msf{d}_A^\fp(\lbrack T' \rbrack)$ have a common constituent. It follows from part \ref{dec_block_compat:b} that $S'$ and $T'$ lie in the same $\fp$-family of $A^K$. Since $S'$ is in the same $A^K$-family as $S$, it is also in the same $\fp$-family as $S$ because the $\fp$-families are unions of $A^K$-families. Similarly, $T'$ is in the same $\fp$-family as $T$. Hence, $S$ and $T$ lie in the same $\fp$-family.

Conversely, suppose that $S$ and $T$ lie in the same $\fp$-family. We have to show that they lie in the same Brauer $\fp$-family. Let $(S_i)_{i=1}^n$ be a system of representatives of the isomorphism classes of simple $A^K$-modules and let $(U_j)_{j=1}^m$ be a system of representatives of the isomorphism classes of simple $A(\fp)$-modules. Let $\sQ \dopgleich (Q_i)_{i=1}^n$ with $Q_i$ being the projective cover of $S_i$, and let $\sP \dopgleich (P_j)_{j=1}^m$ with $P_j$ being the projective cover of $U_j$. Let $\msf{C}_{A(\fp)}$ be the matrix of the Cartan map $\msf{c}_{A(\fp)}$ with respect to the chosen bases, and similarly let $\msf{C}_{A^K}$ be the matrix of $\msf{c}_{A^K}$. Furthermore, let $\msf{D}_A^\fp$ be the matrix of $\msf{d}_A^\fp$ with respect to the chosen bases. Since $\msf{C}_{A(\fp)} = \msf{D}_A^\fp \msf{C}_{A^K} (\msf{D}_A^\fp)^{\rT}$ by Brauer reciprocity, Theorem \ref{brauer_rec}, we have
\begin{equation} \label{brauer_p_fam_thm_equ}
(\msf{C}_{A(\fp)})_{p,q} = (\msf{D}_A^\fp \msf{C}_{A^K} (\msf{D}_A^\fp)^{\rT})_{p,q} = \sum_{k,l=1}^n (\msf{D}_A^\fp)_{p,k} (\msf{C}_{A^K})_{k,l} (\msf{D}_A^{\fp})_{q,l} 
\end{equation}
for all $p,q$. 
Let $U$ be a constituent of $\msf{d}_A^\fp(\lbrack S \rbrack)$ and let $V$ be a constituent of $\msf{d}_A^\fp(\lbrack T \rbrack)$. Since $S$ and $T$ lie in the same $\fp$-family of $A^K$, both $\msf{d}_A^\fp(\lbrack S \rbrack)$ and $\msf{d}_A^\fp(\lbrack T \rbrack)$ lie in the same block of $A(\fp)$ by \ref{dec_block_compat:b}, and therefore $U$ and $V$ lie in the same family of $A(\fp)$. As the families of $A(\fp)$ are equal to the $\sP$-families of $A(\fp)$ by \S\ref{notations}, there exist functions $f: \lbrack 1,r \rbrack \rarr \lbrack 1,m \rbrack$, $g: \lbrack 1,r-1 \rbrack \rarr \lbrack 1,m \rbrack$ with the following properties: $U_{f(1)} = U$, $U_{f(r)} = V$, and for any $j \in \lbrack 1,r-1 \rbrack$ both $P_{f(j)}$ and $P_{f(j+1)}$ have $U_{g(j)}$ as a constituent. We can visualize the situation as follows:
\[
\begin{tikzcd}[column sep=tiny]
P_{f(j)} & & P_{f(j+1)}  \\
U_{f(j)} \arrow{u} & U_{g(j)} \arrow{ul} \arrow{ur} & U_{f(j+1)}  \arrow{u}
\end{tikzcd}
\]
where an arrow $U \rarr P$ signifies that $U$ is a constituent of $P$. 
For any $j \in \lbrack 1,r-1 \rbrack$ we have $(\msf{C}_{A(\fp)})_{g(j),f(j)} \neq 0$ and so it follows from (\ref{brauer_p_fam_thm_equ}) that there are indices $k(j)$ and $l(j)$ such that
\[
(\msf{D}_A^\fp)_{g(j),k(j)} \neq 0\;, \quad (\msf{C}_{A^K})_{k(j),l(j)} \neq 0 \;, \quad (\msf{D}_A^\fp)_{f(j),l(j)} \neq 0 \;.
\] 
Similarly, since $(\msf{C}_{A(\fp)})_{g(j),f(j+1)} \neq 0$, there exist indices $k'(j)$ and $l'(j)$ such that
\[
(\msf{D}_A^\fp)_{g(j),k'(j)} \neq 0\;, \quad (\msf{C}_{A^K})_{k'(j),l'(j)} \neq 0 \;, \quad (\msf{D}_A^\fp)_{f(j+1),l'(j)} \neq 0 \;.
\]
This can be visualized as follows:
\[
\begin{tikzcd}[column sep=tiny]
\msf{d}_A^\fp(\lbrack S_{l(j)} \rbrack) \arrow[dashed,-]{rr} & & \msf{d}_A^\fp(\lbrack S_{k(j)} \rbrack) & & \msf{d}_A^\fp(\lbrack S_{k'(j)} \rbrack) \arrow[dashed,-]{rr} & & \msf{d}_A^\fp(\lbrack S_{l'(j)} \rbrack) \\
U_{f(j)} \arrow{u} & & & U_{g(j)} \arrow{ul} \arrow{ur} & & & U_{f(j+1)} \arrow{u}  
\end{tikzcd}
\]
Here, the dashed edges in the upper row signify that the respective simple $A^K$-modules lie in the same $A^K$-family. Since $U_{f(1)} = U$ and $U_{f(r)} = V$, this shows that $S$ and $T$ lie in the same Brauer $\fp$-family of $A^K$. \qedhere
\end{enum_proof}

\end{proof}

\section{Semicontinuity of blocks in case of a non-split generic fiber} \label{non_split_semicont}

Let $A$ be a finite flat algebra over an integral domain $R$ with fraction field $K$. We have a map
\begin{equation} \label{block_number_prime_map}
\# \Bl_A' : \Spec(R) \to \bbN \;, \quad \fp \mapsto \# \Bl(A(\fp)) \;,
\end{equation}
i.e., $\#\Bl_A'(\fp)$ is the number of blocks of the \textit{specialization} $A(\fp)$. Recall that in \eqref{block_number_map} we considered the map $\#\Bl_A$ with $\#\Bl_A(\fp) = \#\Bl(A_\fp)$ being the number of blocks of the \textit{localization} $A_\fp$. In case $R$ is normal and $A^K$ splits, we know from Proposition \ref{unibranched_block_bijective} that $\#\Bl_A' = \#\Bl_A$. In particular, the map $\# \Bl_A'$ is lower semicontinuous and thus defines a stratification of $\Spec(R)$, the block number stratification, see \S\ref{sec_block_num_strat}.

In case $A^K$ does not split, it still makes perfect sense to consider the map \eqref{block_number_prime_map} and ask if it is lower semicontinuous so that we have a stratification of $\Spec(R)$ by the number of blocks of specializations. But since we do not have the connection from Proposition \ref{unibranched_block_bijective} between blocks of localizations of and blocks of specializations any more, we cannot directly apply the results from \S\ref{blocks_of_localizations}. 

In this section, we will establish a setting where the map $\# \Bl_A'$ is still lower semicontinuous without assuming that the generic fiber $A^K$ splits, the main result being Corollary \ref{finite_type_alg_closed_setting}. To achieve this, however, we have to restrict this map to a subset of $\Spec(R)$. As we will see, in general it is not possible to have lower semicontinuity on all of $\Spec(R)$. \\

First of all, because of the difference between blocks of localizations and blocks of specializations, we introduce the following sets:
\begin{equation}
\beta(A) \dopgleich \msf{max}\lbrace \#\msf{Bl}(A(\fp)) \mid \fp \in \Spec(R) \rbrace \;,
\end{equation}
\begin{equation}
\msf{BlEx}'(A) \dopgleich \#\msf{Bl}_A'^{-1}(\leq \! \beta(A)-1) \;,
\end{equation}
\begin{equation}
\msf{BlGen}'(A) \dopgleich \Spec(R) \setminus \msf{BlEx}'(A) \;.
\end{equation}
Note that if $R$ is normal and $A^K$ splits, then $\beta(A) = \#\Bl(A^K)$, so $\msf{BlEx}'(A) = \msf{BlEx}(A)$ and $\msf{BlGen}'(A) = \msf{BlGen}(A)$ as defined in \eqref{blex_def} and \eqref{blgen_def}. \\

Now, assume that $R'$ is an integral extension of $R$ which is also an integral domain. Let $K'$ be the fraction field of $R'$ and let $\psi:\Spec(R') \twoheadrightarrow \Spec(R)$ be the morphism induced by $R \subs R'$. The scalar extension $A' \dopgleich R' \otimes_{R} A$ is again a finitely presented flat $R'$-algebra (using Remark \ref{fp_projective_remark}). For any $\fp \in \Spec(R)$ and any $\fp' \in \Spec(R')$ lying over $\fp$ we have a diagram

\begin{equation} \label{prop_P1_morphisms}
\begin{tikzcd}
 & A_{\fp'}' \arrow[twoheadrightarrow]{d} \\
A(\fp) = A_\fp/\fp_\fp A_\fp \arrow[hookrightarrow]{r} & A'(\fp') = A'_{\fp'}/{\fp'}_{\fp'} A'_{\fp'} 
\end{tikzcd}
\end{equation}
and it then follows from (\ref{get_more_blocks_equation}) that
\begin{equation} \label{main_thm_proof_block_lesseq}
\#\msf{Bl}(A(\fp)) \leq \#\msf{Bl}(A'(\fp')) \geq \#\msf{Bl}(A'_{\fp'}) \;.
\end{equation}
Let $X$ be a set contained in
\begin{equation}
X_{R'}(A) \dopgleich \lbrace \fp \in \Spec(R) \mid \#\msf{Bl}(A(\fp)) \!=\! \#\msf{Bl}(A'(\fp')) \!=\! \#\msf{Bl}(A'_{\fp'}) \tn{ for all } \fp' \in \psi^{-1}(\fp) \rbrace .
\end{equation}
We have seen in Corollary \ref{normal_split_corollary} that in case $R$ is normal and $A^K$ splits we can choose $R=R'$ and have $X = \Spec(R)$. In general $X$ will be a proper subset of $\Spec(R)$ and we have to choose $R'$ appropriately to enlarge it a bit more. Let us first concentrate on what we can say when restricting to $X$. We introduce the following restricted versions of our invariants:
\begin{equation}
\#\Bl_{A,X}' \dopgleich \#\Bl_A' |_X : X \to \bbN \;,
\end{equation}
\begin{equation} \label{bl_leqn_blowdown}
\# \Bl_{A,X}'^{-1}(\leq\! n) \dopgleich \#\Bl_{A}'^{-1}(\leq \! n ) \cap X = \psi(\# \Bl_{A'}^{-1}(\leq \! n) )\cap X \;,
\end{equation}
\begin{equation}
\# \Bl_{A,X}'^{-1}( n) \dopgleich \#\Bl_{A}'^{-1}( n ) \cap X = \psi(\# \Bl_{A'}^{-1}( n) )\cap X \;,
\end{equation}
\begin{equation}
\beta_X(A) \dopgleich \msf{max}\lbrace \#\Bl(A(\fp)) \mid \fp \in X \rbrace \;,
\end{equation}
\begin{equation}
\msf{BlEx}'_X(A) \dopgleich \Bl_{A,X}'^{-1}(\leq \! \beta_X(A)-1) \;,
\end{equation}
\begin{equation}
\msf{BlGen}_X'(A) \dopgleich X \setminus \msf{BlEx}'_X(A)  \;.
\end{equation}

\begin{cor} \label{bln_spec_strat}
The map $X \rarr \bbN$, $\fp \mapsto \#\msf{Bl}(A(\fp))$, is lower semicontinuous on $X$, so  $X = \coprod_{n \in \bbN} \msf{Bl}_{A,X}'^{-1}(n)$ is a stratification of $X$ into locally closed subsets. Moreover, 
\begin{equation} \label{beta_invariant_estimate}
\beta_X(A) \leq \#\Bl(A^{K'}) \;.
\end{equation}
\end{cor}

\begin{proof}
Since $\psi$ is a closed morphism and $\#\Bl_{A'}^{-1}(\leq \! n)$ is closed in $\Spec(R')$ by (\ref{bl_leqn_closed}), it follows that $\psi(\#\Bl_{A'}^{-1}(\leq \! n)))$ is closed in $\Spec(R)$, hence $\# \Bl_{A,X}'^{-1}(\leq \! n )$ is closed in $X$ by (\ref{bl_leqn_blowdown}). Since $\#\Bl_{A,X}'^{-1}(n)  = \#\Bl_{A,X}'^{-1}(\leq \! n)  \setminus  \#\Bl_{A,X}'^{-1}(\leq \! n-1)$, it is clear that $\#\Bl_{A,X}'^{-1}(n)$ is locally closed in $X$. We have seen in (\ref{bln_closure}) that $\ol{\#\Bl_{A'}^{-1}(n)} \subs \bigcup_{m \leq n} \Bl_{A'}^{-1}(m)$. Hence, since $\psi$ is closed, we obtain
\[
\ol{\# \Bl_{A,X}'^{-1}(n)} = \psi(\ol{\#\Bl_{A'}^{-1}(n)}) \cap X \subs \bigcup_{m \leq n} \psi(\# \Bl_{A'}^{-1}(m)) \cap X = \bigcup_{m \leq n} \#\Bl_{A,X}'^{-1}(m) \;. 
\]
\end{proof}

Note that in (\ref{beta_invariant_estimate}) we could only bound $\beta_X(A)$ above by $\#\Bl(A^{K'})$, and not by $\# \Bl(A^K)$. In fact, we will see in Example \ref{ken_example} that we may indeed have $\beta_X(A) > \#\Bl(A^K)$ in general. This is an important difference to blocks of localizations where we always have the maximal number of blocks in the generic point. 

In the following lemma we describe a situation where we have $\beta_X(A) = \#\Bl(A^{K'})$. We recall that $X$ is called \textit{very dense} if the embedding $X \hookrightarrow \Spec(R)$ is a \textit{quasi-homeomorphism}, i.e., the map $Z \mapsto Z \cap X$ is a bijection between the closed (equivalently, open) subsets of the two spaces. This notion was introduced by Grothendieck \cite[\S10]{Grothendieck:EGA-4-3}.

\begin{lem} \label{x_very_dense_beta}
Suppose that $X$ is \textnormal{very dense} in $\Spec(R)$, that $R$ is noetherian, and that $\psi$ is finite. Then $\beta_X(A) = \#\Bl(A^{K'})$, thus $\msf{BlEx}_X'(A) = \psi(\msf{BlEx}(A')) \cap X$. If moreover $R'$ is normal and $R$ is universally catenary, then $\msf{BlEx}_X'(A)$ is a reduced Weil divisor in $X$.
\end{lem}

\begin{proof}
The assumptions imply that $R'$ is noetherian, too. We know from Theorem \ref{block_stratification_thm} that $\msf{BlGen}(A')$ is a non-empty open subset of $\Spec(R')$. In particular, it is constructible. Since $\Spec(R)$ is quasi-compact, the morphism $\psi$ is quasi-compact by \cite[Remark 10.2.(1)]{GorWed10-Algebraic-geomet}. It thus follows from Chevalley's constructibility theorem, see \cite[Corollary 10.71]{GorWed10-Algebraic-geomet}, that $\psi(\msf{BlGen}(A'))$ is constructible in $\Spec(R)$. Since $X$ is very dense in $\Spec(R)$, we conclude that $\psi(\msf{BlGen}(A')) \cap X \neq \emptyset$ by \cite[Proposition 10.1.2]{Grothendieck:EGA-4-3}. Hence, there is $\fp \in X$ and $\fp' \in \msf{BlGen}(A')$ with $\psi(\fp') = \fp$. But then we have $\#\Bl(A(\fp)) = \#\Bl(A'(\fp')) = \#\Bl(A^{K'})$, so $\beta_X(A) = \#\Bl(A^{K'})$. Now, assume that $R'$ is normal and $R$ is universally catenary. We know that $\msf{BlEx}(A')$ is either empty or pure of codimension one in $\Spec(R')$ by Corollary \ref{blex_is_weil_div}. In \cite[Theorem B.5.1]{Huneke-Swanson-Integral-Closure} it is shown that the extension $R \subs R'$ satisfies the dimension formula, hence $\psi(\msf{BlEx}(A'))$ is either empty or pure of codimension one. Since $X$ is very dense in $\Spec(R)$, the same is also true for $X \cap  \psi(\msf{BlEx}(A')) =\msf{BlEx}_X'(A)$.
\end{proof}

\begin{cor} \label{finite_type_alg_closed_setting}
Suppose that $R$ is a finite type algebra over an algebraically closed field. Let $X$ be the set of closed points of $\Spec(R)$. Then the map $X \rarr \bbN$, $\fm \mapsto \#\msf{Bl}(A(\fm))$, is lower semicontinuous and so $X = \coprod_{n \in \bbN} \# \Bl_{A,X}'^{-1}(n)$ is a stratification of $X$. Moreover, $\beta_X(A) = \#\Bl(A^{\ol{K}})$, where $\ol{K}$ is an algebraic closure of $K$. If $R$ is also universally catenary, then $\msf{BlEx}_X'(A)$ is a reduced Weil divisor in $X$. 
\end{cor}

\begin{proof}
Let $K'$ be a finite extension of $K$ such that $A^{K'}$ splits (this is always possible, see \cite[Proposition 7.13]{CR-Methods-1}) and let $R'$ be the integral closure of $R$ in $K'$. Now, $\#\Bl(A'(\fp')) = \#\Bl(A'_{\fp'})$ for all $\fp' \in \Spec(R)$ by Proposition \ref{unibranched_block_bijective}. Since $R$ is a finite type algebra over an algebraically closed field $k$, the residue field in a closed point $\fm$ of $\Spec(R)$ is just $k$. Hence, the specialization $A(\fm)$ is a finite-dimensional algebra over an algebraically closed field, thus splits and we therefore have $\#\msf{Bl}(A(\fm)) = \#\msf{Bl}(A'(\fm'))$ for any $\fm' \in \psi^{-1}(\fm)$ by Lemma \ref{split_center_lemma}. Hence, $X \subs X_{R'}(A)$. The claim about semicontinuity and the stratification thus follows from Corollary \ref{bln_spec_strat}. It is shown in \cite[Proposition 3.35]{GorWed10-Algebraic-geomet} that $X$ is very dense in $\Spec(R)$. Since $R$ is a finite type algebra over a field, it is  japanese, so $\psi$ is a finite morphism. Hence, $\beta_X(A) = \#\msf{Bl}(A^{K'}) = \#\msf{Bl}(A^{\ol{K}})$ by Lemma \ref{x_very_dense_beta}. The claim that $\msf{BlEx}_X'(A)$ is a reduced Weil divisor if $R$ is universally catenary also follows from Lemma \ref{x_very_dense_beta}.
\end{proof}

\begin{ex} \label{ken_example}
The following example due to K. Brown shows that in the setting of Lemma \ref{finite_type_alg_closed_setting} we may indeed have $\beta_X(A) > \#\msf{Bl}(A^K)$ so that the map $\fp \mapsto \#\msf{Bl}(A(\fp))$ will not be lower semicontinuous on the whole of $\Spec(R)$. Let $k$ be an algebraically closed field of characteristic zero, let $X$ be an indeterminate over $k$, let $R \dopgleich k \lbrack X^n \rbrack$ for some $n > 1$, and let $A \dopgleich k \lbrack X \rbrack$. Let $C_n$ be the cyclic group of order $n$. We fix a generator of $C_n$ and let it act on $X$ by multiplication with a primitive $n$-th root of unity. Then $R = k \lbrack X \rbrack^{C_n}$, so $A$ is free of rank $n$ over $R$. Moreover, $\Frac(A) = k(X)$ is a Galois extension of degree $n$ of $K \dopgleich \Frac(R)$ by \cite[Proposition 1.1.1]{Ben-Polynomial-invariants}, so in particular  $K \neq k(X)$ since $n>1$. By \cite[Ex. 6R]{Goodearl-Warfield} we have 
\[
A^K = A \otimes_R K = A\lbrack (R \setminus \lbrace 0 \rbrace)^{-1}\rbrack = \Frac(A) = k(X) \;,
\]
so the $K$-algebra $A^K = \msf{Z}(A^K)$ is \textit{not} split (and thus also not block-split by Lemma \ref{split_center_lemma}). It is clear that 
\begin{equation}
\#\msf{Bl}(A^K) = 1 \;.
\end{equation}
Now, let $\fm \dopgleich (X^n-1) \in \msf{Max}(R)$. Then $\msf{k}(\fp) = k$ and since $k$ is algebraically closed, we have $A(\fm) = A/\fm A \simeq k^n$ as $k$-algebras. In particular, 
\begin{equation}
\#\msf{Bl}(A(\fm)) = n > 1 = \#\msf{Bl}(A^K) \;.
\end{equation}
\end{ex}

We close with a setting where our base ring is not necessarily normal but we still get a global result on $\Spec(R)$.%

\begin{lem} \label{split_fibers_results}
Suppose that $A$ has split fibers, i.e., $A(\fp)$ splits for all $\fp \in \Spec(R)$. Then the map $\Spec(R) \rarr \bbN$, $\fp \mapsto \#\msf{Bl}(A(\fp))$, is lower semicontinuous and so $\Spec(R) = \coprod_{n \in \bbN} \msf{Bl}_A'^{-1}(n)$ is a partition into locally closed subsets. Moreover, $\beta(A) = \#\msf{Bl}(A^K)$. If $R$ is also  universally catenary, japanese, and noetherian, then $\msf{BlEx}'(A)$ is a reduced Weil divisor in $\Spec(R)$.%
\end{lem}

\begin{proof}
Let $R'$ be the integral closure of $R$ in $K$. Then $\#\msf{Bl}(A'(\fp')) = \#\msf{Bl}(A'_{\fp'})$ for all $\fp' \in \Spec(R')$ by Proposition \ref{unibranched_block_bijective}. Since $A(\fp)$ splits, we moreover have $\#\msf{Bl}(A(\fp)) = \#\msf{Bl}(A'(\fp'))$ for all $\fp \in \Spec(R)$ $\fp' \in \psi^{-1}(\fp)$ by Lemma \ref{split_center_lemma}. Hence, $X_{R'}(A) = \Spec(R)$. The claim about semicontinuity and the partition follows from Corollary \ref{bln_spec_strat}. Now, assume that $R$ is universally catenary, japanese, and noetherian. Since $R$ is japanese, it follows by definition that $\psi$ is finite. The claim about $\msf{BlEx}'(A)$ being a reduced Weil divisor now follows from \ref{x_very_dense_beta}. %
\end{proof}

\appendix

\section{More on base change of blocks} \label{appendix_blocks_base_change}

In this appendix we collect several facts about base change of blocks. Some results here should also be of independent interest.

\subsection{Block compatibility of scalar extension of modules}
Recall the decomposition of the module category of a ring $A$ relative to a decomposition in of $1 \in A$ into pairwise orthogonal central idempotents described in \S\ref{notations}. We have the following compatibility.

\begin{lem} \label{block_extension_compatibility}
Let $\phi:R \rarr S$ be a morphism of commutative rings and let $A$ be an $R$-algebra. Suppose that $\phi_A$ is central idempotent stable and let $V$ be a non-zero $A$-module. In any of the following cases the $A$-module $V$ belongs to the block $c_i$ if and only if the $A^S$-module $V^S$ belongs to the $\phi$-block $\phi_A(c_i)$:
\begin{enum_thm}
\item $\phi$ is injective and $V$ is $R$-projective.
\item $\phi$ is faithfully flat.
\item $R$ is local or a principal ideal domain and $V$ is $R$-free.
\end{enum_thm}
\end{lem}

\begin{proof}
As $c_jV$ is a direct summand of $V$, it follows that we have a canonical isomorphism between isomorphism $\phi_A^*(c_j V) \simeq \phi_A(c_j)\phi_A^*V $ of $A^S$-modules for all $j$. The claim thus holds if we can show that no non-zero direct summand $V'$ of $V$ is killed by $\phi_A^*$, i.e., $\phi_A^*V' \neq 0$. But this is implied by the assumptions in each case. Namely, in the first two cases it follows from Lemma \ref{phi_injective_lemma} that $\phi_V$ is injective, which implies that $\phi_{V'}$ is also injective, so $\phi_A^*V'$ cannot be zero for non-zero $V'$. In the third case neither $\phi$ nor $\phi_V$ need to be injective, so this needs extra care. First of all, since $V$ is assumed to be $R$-free, the assumptions on $R$ imply that a direct summand $V'$ of $V$, which a priori is only $R$-projective, is already $R$-free, too. In case $R$ is local, this follows from Kaplansky's theorem \cite{Kaplansky} and in case $R$ is a principal ideal domain, this is a standard fact. Now, if $V$ is $R$-free with basis $( v_\lambda )_{\lambda \in \Lambda}$, then it is a standard fact (see \cite[II, \S5.1, Proposition 4]{Bou-Algebra-1-3}) that $\phi_A^*V$ is $S$-free with basis $(\phi_V(v_\lambda))_{\lambda \in \Lambda}$. This shows that $\phi_A^*V \neq 0$ for any non-zero $R$-free $A$-module $V$. This applied to direct summands of $V$, which are $R$-free as shown, proves the claim.
\end{proof}

\subsection{Field extensions} \label{blocks_field_extensions}

Throughout this paragraph let $A$ be a finite-dimensional algebra over a field $K$. From (\ref{get_more_blocks_equation}) we know that $\#\msf{Bl}(A) \leq \#\msf{Bl}(A^L)$ for any extension field $L$ of $K$.

\begin{defn} \label{block_split_def}
We say that $A$ is \word{block-split} if $\#\msf{Bl}(A) = \#\msf{Bl}(A^L)$ for any extension field $L$ of $K$.
\end{defn}

Our aim is to show the following lemma.

\begin{lem} \label{split_center_lemma}
If $\msf{Z}(A)$ is a split $K$-algebra (e.g., if $A$ itself splits), then $A$ is block-split. The converse holds if $K$ is perfect.
\end{lem}

The first assertion of the lemma is essentially obvious since $\msf{Z}(A)$ is semiperfect and therefore
\begin{equation} \label{bl_irr_equation}
\#\msf{Bl}(A) = \#\msf{Bl}(\msf{Z}(A)) = \msf{rk}_\bbZ \ \msf{K}_0(\msf{Z}(A)) = \#\msf{rk}_\bbZ \ \msf{G}_0(\msf{Z}(A)) = \#\Irr \msf{Z}(A) \;,
\end{equation}
where the second equality follows from the fact that idempotents in a commutative ring are isomorphic if and only if they are equal, see \cite[Ex. 22.2]{Lam-First-Course-91}. The same equalities of course also hold for $\msf{Z}(A)^L = \msf{Z}(A^L)$, where $L$ is an extension field of $K$. Hence, if $\msf{Z}(A)$ is split, then $A$ is block-split. If $A$ itself is split, it is a standard fact that its center splits, so $A$ is block-split. \\

We will prove the converse (assuming that $K$ is perfect) from a more general point of view as the results might be of independent interest and we re-use some of them in the last section. First of all, the field extension $K \subs L$ induces natural group morphisms
\begin{equation}
\msf{d}_A^L:\msf{G}_0(A) \rarr \msf{G}_0(A^L) \quad \tn{and} \quad 
\msf{e}_A^L:\msf{K}_0(A) \rarr \msf{K}_0(A^L) \;.
\end{equation}
Without any assumptions on the field $K$ we have the following property.

\begin{lem} \label{d_e_injective}
The morphisms $\msf{d}_A^L$ and $\msf{e}_A^L$ are injective.
\end{lem}

\begin{proof}
Let $(S_i)_{i \in I}$ be a system of representatives of the isomorphism classes of simple $A$-modules. For each $i$ let $(T_{ij})_{j \in J_i}$ be a system of representatives of the isomorphism classes of simple $A^L$-modules which occur as constituents of $S_i^L$. Then by \cite[Proposition 7.13]{Lam-First-Course-91} the set $(T_{ij})_{i \in I, j \in J_i}$ is a system of representatives of the isomorphism classes of simple $A^L$-modules. Hence, the matrix $\msf{D}_A^L$ of $\msf{d}_A^L$ in bases given by the isomorphism classes of simple modules is in column-echelon form, has no zero columns, and no zero rows. In particular, $\msf{d}_A^L$ is injective.

For each $i \in I$ let $P_i$ be the projective cover of $S_i$ and for each $j \in J_i$ let $Q_{ij}$ be the projective cover of $T_{ij}$. By the above, $(Q_{ij})_{i \in I, j \in J_i}$ is a system of representatives of the isomorphism classes of projective indecomposable $A^L$-modules. We claim that in the direct sum decomposition of the finitely generated projective $A^L$-module $P_i^L$ into projective indecomposable $A^L$-modules only the $Q_{ij}$ with $j \in J_i$ occur. With the same argument as above, this implies that $\msf{e}_A^L$ is injective. So, let us write $P_i^L = \bigoplus_{\lambda \in \Lambda} U_\lambda$ for (not necessarily non-isomorphic) projective indecomposable $A^L$-modules $U_\lambda$. The $U_\lambda$ are the up to isomorphism unique projective indecomposable $A^L$-modules occurring as direct summands of $P_i^L$. As the radical is additive by \cite[Proposition 24.6(ii)]{Lam-First-Course-91}, we have $\Rad(P_i^L) = \bigoplus_{\lambda \in \Lambda} \Rad(U_\lambda)$, so 
\begin{equation}
S_i^L = (P_i/\Rad(P_i))^L = P_i^L/\Rad(P_i)^L = \bigoplus_{\lambda \in \Lambda} U_\lambda/(\Rad(P_i)^L \cap U_\lambda) \;.
\end{equation}
Moreover, we have $\Rad(P_i)^L \subs \Rad(P_i^L)$. This follows from the fact that $\Rad(A)^L \subs \Rad(A^L)$ by \cite[Theorem 5.14]{Lam-First-Course-91} and the fact that $\Rad(P_i) = \Rad(A) P_i$ and $\Rad(P_i^L) = \Rad(A^L) P_i^L$ by \cite[Theorem 24.7]{Lam-First-Course-91} since $P_i$ and $P_i^L$ are projective. For each $\lambda \in \Lambda$ the radical of $U_\lambda$ is a proper submodule of $U_\lambda$ and therefore
\[
\Rad(P_i)^L \cap U_\lambda \subs \Rad(P_i^L) \cap U_\lambda = \Rad(U_\lambda) \subsetneq U_\lambda \;. 
\]
Hence, the head of $U_\lambda$ is a constituent of $U_\lambda/(\Rad(P_i^L) \cap U_\lambda)$, and since all constituents of the latter are constituents of $S_i^L$, we must have $\msf{Hd}(U_\lambda) \simeq S_{ij_{\lambda}}$ for some $j_\lambda \in J_i$ by the above. This implies that $U_\lambda = Q_{ij_\lambda}$, thus proving the claim.
\end{proof}

\begin{lem} \label{d_iso_iff_e_iso}
The following holds:
\begin{enum_thm}
\item The morphism $\msf{d}_A^L$ is an isomorphism if and only if it induces a bijection between isomorphism classes of simple modules. Similarly, the morphism $\msf{e}_A^L$ is an isomorphism if and only if it induces a bijection between isomorphism classes of projective indecomposable modules.
\item \label{d_iso_iff_e_iso:reflect} If $\msf{d}_A^L$ is an isomorphism, so is $\msf{e}_A^L$. The converse holds if $K$ is perfect.
\end{enum_thm}
\end{lem}

For the proof of Lemma \ref{d_iso_iff_e_iso} we will need the following well-known elementary lemma that is also used in the last section. Recall from (\ref{intertwining_form_def}) the intertwining form $\langle \cdot,\cdot \rangle_A$ of $A$.

\begin{lem} \label{intertwining_nondeg}
Let $P$ be a projective indecomposable $A$-module and let $V$ be a finitely generated $A$-module. Then
\begin{equation}
\langle \lbrack P \rbrack, \lbrack V \rbrack \rangle_A = \lbrack V:\msf{Hd}(P) \rbrack \cdot \dim_K \End_A(\msf{Hd}(P)) \;,
\end{equation}
where $\msf{Hd}(P) = P/\Rad(P)$ is the head of $P$. In particular, $\langle \cdot, \cdot \rangle_A$ is non-degenerate.
\end{lem}

\begin{proof}
We first consider the case $V=\msf{Hd}(P)$. Let $f \in \Hom_A(P,\Hd(P))$ be non-zero. Since $\Hd(P)$ is simple, this morphism is already surjective and thus induces an isomorphism $P/\Ker(f) \cong \Hd(P)$. But as $\Rad(P)$ is the unique maximal submodule of $P$, we must have $\Ker(f) = \Rad(P)$ and thus get an induced morphism $\Hd(P) \rarr \Hd(P)$. This yields a $K$-linear morphism $\Phi: \Hom_A(P,\Hd(P)) \rarr \End_A(\Hd(P))$. On the other hand, if $f \in \End_A(\Hd(P))$, then composing it with the quotient morphism $P \rarr P/\Rad(P) = \Hd(P)$ yields a morphism $P \rarr \Hd(P)$. In this way we also get a $K$-linear morphism $\Psi:\End_A(\Hd(P)) \rarr  \Hom_A(P,\Hd(P))$. By construction, $\Phi$ and $\Psi$ are pairwise inverse, hence $\langle \lbrack P \rbrack, \lbrack \msf{Hd}(P) \rbrack \rangle_A = \dim_K \Hom_A(P,\msf{Hd}(P)) = \dim_K \End_A(\Hd(P))$ as claimed.

Now, suppose that $V$ is a simple $A$-module not isomorphic to $\msf{Hd}(P)$. We can write $P = Ae$ for some primitive idempotent $e \in A$. Since $A$ is artinian, $e$ is already local and now it follows from \cite[21.19]{Lam-First-Course-91} that $\Hom_A(Ae, V)$ is non-zero if and only if $V$ has a constituent isomorphic to $\Hd(Ae)$. This is not true by assumption, and therefore $\Hom_A(P,V) = 0$, so $\langle \lbrack P \rbrack, \lbrack V \rbrack \rangle_A = 0$.

Finally, for $V$ general we have $
\lbrack V \rbrack = \sum_{S \in \Irr A} \lbrack V : S \rbrack \lbrack S \rbrack$ in $\msf{G}_0(A)$. By the above we get 
\begin{align*}
\langle \lbrack P \rbrack, \lbrack V \rbrack \rangle_A &= \sum_{S \in \Irr A} \lbrack V:S \rbrack \langle \lbrack P \rbrack, \lbrack S \rbrack \rangle_A = \lbrack V:\Hd(P) \rbrack \langle \lbrack P \rbrack, \lbrack \Hd(P) \rbrack \rangle_A \\ &= \lbrack V : \Hd(P) \rbrack \cdot \dim_K \End_A(\Hd(P)) 
\end{align*}
as claimed. It follows that the Gram matrix $\mscr{G}$ of $\langle \cdot,\cdot \rangle$ with respect to the basis $(\msf{P}(S))_{S \in \Irr A}$ of $\msf{K}_0(A)$ and the basis $(S)_{S \in \Irr A}$ of $\msf{G}_0(A)$ is diagonal with positive diagonal entries. The determinant of $\mscr{G}$ is thus a non-zero divisor on $\bbZ$ and since both $\msf{K}_0(A)$ and $\msf{G}_0(A)$ are $\bbZ$-free of the same finite dimension, it follows that $\langle \cdot,\cdot \rangle_A$ is non-degenerate, see \cite[Satz 70.5]{Scheja-Storch-Alg-2}. 
\end{proof}

\begin{proof}[Proof of Lemma \ref{d_iso_iff_e_iso}]
We use the same notations as in the proof of Lemma \ref{d_e_injective}. Since $_AA$ is a projective $A$-module, there is a decomposition $_AA = \bigoplus_{i \in I} P_i^{r_i}$ for some $r_i \in \bbN$.  
Using Lemma \ref{intertwining_nondeg} we see that
\begin{align*}
\dim_K \Hd(P_j) & = \langle \lbrack _AA \rbrack, \lbrack \Hd(P_j) \rbrack \rangle_A = \sum_{i \in I} r_i \langle \lbrack P_i \rbrack, \lbrack \Hd(P_j) \rbrack \rangle_A = r_j \langle \lbrack P_j \rbrack, \lbrack \Hd(P_j) \rbrack \rangle_A \\
& = r_j \dim_K \End_A(\Hd(P_i)) \;.
\end{align*}
Hence, $r_i = \frac{n_i}{m_i}$, where $n_i \dopgleich \dim_K S_i$ and $m_i \dopgleich \dim_K \End_A(S_i)$. In particular, 
\begin{equation} \label{d_iso_iff_e_iso:equ}
\dim_K A = \sum_{i \in I} \frac{n_i}{m_i} \dim_K P_i \;.
\end{equation}

Now, suppose that $\msf{d}_A^L$ is an isomorphism. Then clearly $\#\Irr A = \#\Irr A^L$. The properties of the matrix $\msf{D}_A^L$ of the morphism $\msf{d}_A^L$ derived in the proof of Lemma \ref{d_e_injective} immediately imply that $\msf{D}_A^L$ is diagonal. Since it is invertible with natural numbers on the diagonal, it must already be the identity matrix, i.e., $\msf{d}_A^L$ induces a bijection between the isomorphism classes of simple modules. In particular, $(S_i^L)_{i \in I}$ is a system of representatives of the isomorphism classes of simple $A^L$-modules. The properties of the matrix $\msf{E}_A^L$ of $\msf{e}_A^L$ derived in the proof of Lemma \ref{d_e_injective} now imply that we must have $P_i^L \simeq Q_i^{s_i}$ for some $s_i \in \bbN$. We argue that $s_i = 1$. This shows that $\msf{e}_A^L$ is an isomorphism inducing a bijection between the isomorphism classes of projective indecomposable modules. In the same way we deduced equation (\ref{d_iso_iff_e_iso:equ}) we now get
\begin{equation} \label{d_iso_iff_e_iso:equ2}
\dim_K A = \dim_L A^L = \sum_{i \in I} \frac{n_i'}{m_i'} \dim_L Q_i\end{equation}
with
\[
n_i' = \dim_L \Hd(Q_i) = \dim_L S_i^L = \dim_K S_i = n_i
\]
and
\[
m_i' = \dim_L \End_{A^L}(\Hd(Q_i)) = \dim_L \End_{A^L}(S_i^L) = \dim_K \End_K(S_i) = m_i \;,
\]
using the fact that $L \otimes_K \End_A(S_i) \simeq \End_{A^L}(S_i^L)$, see Lemma \cite[Theorem 2.38]{Rei-Maximal-orders-0}. Since $\dim_L Q_i \leq \dim_L P_i^L = \dim_K P_i$, equations (\ref{d_iso_iff_e_iso:equ}) and (\ref{d_iso_iff_e_iso:equ2}) imply that $\dim_L Q_i = \dim_K P_i$, so $Q_i = P_i^L$.

Conversely, suppose that $\msf{e}_A^L$ is an isomorphism. With the  properties of the matrix $\msf{E}_A^L$ of $\msf{e}_A^L$ established in the proof of Lemma \ref{d_e_injective} we see similarly as above that $\msf{e}_A^L$ already induces a bijection between the projective indecomposable modules. In particular, $P_i^L \simeq Q_i$. Due to the properties of the matrix $\msf{D}_A^L$ of $\msf{d}_A^L$ established in the proof of Lemma \ref{d_e_injective} the only constituent of $S_i^L$ is $T_i$. Since $P_i$ is the projective cover of $P_i$, we have a surjective morphism $\varphi:P_i \twoheadrightarrow S_i$ with $\Ker(\varphi) = \Rad(P_i)$. Scalar extension induces a surjective morphism $\varphi^L:P_i^L \twoheadrightarrow S_i^L$ with $\Ker(\varphi^L) = \Ker(\varphi)^L = \Rad(P_i)^L \subs \Rad(P_i^L)$. It thus follows from \cite[Corollary 6.25(i)]{CR-Methods-1} that $P_i^L$ is the projective cover of $S_i^L$. Now, we assume that $K$ is perfect. Then by \cite[Theorem 7.5]{CR-Methods-1} all simple $A$-modules are separable, so $S_i^L = T_i^{s_i}$ for some $s_i$. Since projective covers are additive, we get $P_i^L = Q_i^{s_i}$. As $P_i^L = Q_i$, this implies that $s_i=1$, so $S_i^L = T_i$ is simple. Hence, $\msf{d}_A^L$ induces a bijection between the isomorphism classes of simple modules.
\end{proof}

\begin{rem}
With the same arguments as in the proof of Lemma \ref{d_iso_iff_e_iso} we can show that the converse in Lemma \ref{d_iso_iff_e_iso}\ref{d_iso_iff_e_iso:reflect} still holds when we only assume that all simple $A$-modules are separable, i.e., they remain semisimple under field extension. This holds for example when $A$ splits or if $A$ is a group algebra (over any field). We do not know whether it holds more generally.
\end{rem}

\begin{proof}[Proof of Lemma \ref{split_center_lemma}]
Let $Z \dopgleich \msf{Z}(A)$. Suppose that $L$ is an extension field of $K$ with $\# \msf{Bl}(A) = \#\msf{Bl}(A^L)$. By (\ref{bl_irr_equation}) we know that $\#\Irr Z = \#\Irr Z^L$. The arguments in the proof of Lemma \ref{d_e_injective} thus imply that the matrix $\msf{D}_A^L$ of the morphism $\msf{d}_Z^L:\msf{G}_0(Z) \rarr \msf{G}_0(Z^L)$ must be a diagonal matrix. We claim that it is the identity matrix. Since this holds for any $L$, it means that the simple modules of $Z$ remain simple under any field extension, so $Z$ splits. Our assumption implies that $\#\msf{Idem}_p(Z) = \#\msf{Idem}_p(Z^L)$, so every primitive idempotent $e \in Z$ remains primitive in $Z^L$. This shows that $\msf{e}_A^L:\msf{K}_0(Z) \rarr \msf{K}_0(Z^L)$ induces a bijection between projective indecomposable modules. In particular, it is an isomorphism. Now, Lemma \ref{d_iso_iff_e_iso} shows that also $\msf{d}_A^L$ is an isomorphism. Since its matrix $\msf{D}_A^L$ is invertible with natural numbers on the diagonal, it must be the identity.
\end{proof}

\begin{rem}
In the proof of Lemma \ref{split_center_lemma} we have deduced that for a commutative finite-dimensional $K$-algebra $Z$ the condition $\msf{rk}_\bbZ \ \msf{K}_0(Z) = \msf{rk}_\bbZ \ \msf{K}_0(Z^L)$ already implies that $\msf{e}_Z^L$ induces a bijection between projective indecomposable modules. This follows from the fact that idempotents in a commutative ring are isomorphic if and only if they are equal. This is not true for a non-commutative ring $A$. Here, we can have $\msf{rk}_\bbZ \ \msf{K}_0(A) = \msf{rk}_\bbZ \ \msf{K}_0(A^L)$ but still a primitive idempotent $e \in A$ can split into a sum of \textit{isomorphic} orthogonal primitive idempotents of $A^L$. Then the matrix $\msf{E}_A^L$ of $\msf{e}_A^L$ is diagonal but not the identity. %
\end{rem}

Let us record the following additional fact:

\begin{lem} \label{num_blocks_rad_dim}
If $\msf{Z}(A)$ splits, then 
\begin{equation}
\#\msf{Bl}(A) = \dim_K \msf{Z}(A) - \dim_K \Rad(\msf{Z}(A)) = \dim_K \msf{Z}(A) - \dim_K (\msf{Z}(A) \cap \Rad(A) )\;.
\end{equation}
\end{lem}

\begin{proof}
This follows immediately from (\ref{bl_irr_equation}) and the fact that that $\Rad(\msf{Z}(A)) = \msf{Z}(A) \cap \Rad(A)$ since $\msf{Z}(A) \subs A$ is a finite normalizing extension, see \cite[Theorem 1.5]{Lor-Finite-normalizing}.
\end{proof}

\subsection{Reductions}

Now, we consider a situation which in a sense is opposite to the one considered in the last paragraph, namely we consider the quotient morphism $\phi:R \twoheadrightarrow R/\fm \gleichdop S$ for a local commutative ring $R$ with maximal ideal $\fm$ and a finitely generated $R$-algebra $A$. 
By Lemma \ref{idempotent_stable}\ref{idempotent_stable:surj} the morphism $\phi_A:A \twoheadrightarrow A^S \simeq A/\fm A \gleichdop \ol{A}$ is idempotent stable. We say that $\phi$ is \word{idempotent surjective} if for each idempotent $e' \in A^S$ there is an idempotent $e \in A$ with $\phi_A(e) = e'$. We say that $\phi_A$ is \word{primitive idempotent bijective} if it induces a bijection between the \textit{isomorphism classes} of primitive idempotents of $A$ and the isomorphism classes of primitive idempotents of~$A^S$, The question whether $\phi_A$ is idempotent surjective is precisely the question whether idempotents of $\ol{A}$ can be lifted to $A$, and this is a classical topic in ring theory. The following lemma is standard, we omit the proof.

\begin{lem} \label{idempotent_surjective_block_bijective}
If $\phi_A:A \twoheadrightarrow \ol{A}$ is idempotent surjective, it is  primitive idempotent bijective and block bijective.
\end{lem}

The next theorem was proven by M.\ Neunhöffer \cite[Proposition 5.10]{Neunhoeffer-PhD}.

\begin{thm}[M.\ Neunhöffer] \label{neunhoeffer_semiperfect}
The morphism $\phi_A:A \twoheadrightarrow \ol{A}$ is idempotent surjective if and only if $A$ is semiperfect. 
\end{thm}

We recall two standard situations of idempotent surjective reductions.

\begin{lem} \label{semiperfect_std_settings}
In the following two cases the morphism $\phi_A:A \twoheadrightarrow \ol{A}$ is idempotent surjective:
\begin{enum_thm}
\item \label{semiperfect_std_settings:complete} $R$ is noetherian and $\fm$-adically complete.
\item \label{semiperfect_std_settings:hensel} $R$ is henselian.
\end{enum_thm}
\end{lem}

\begin{proof}
For a proof of the first case, see \cite[Proposition 21.34]{Lam-First-Course-91}. For a proof of the second case assuming that $A$ is \textit{commutative}, see \cite[I, \S3, Proposition 2]{Raynaud:Henselian}. To give a proof for non-commutative $A$ let $\ol{e} \in \ol{A}$ be an idempotent. Let $k \dopgleich R/\fm$ and let $\ol{B} \dopgleich k \lbrack \ol{e} \rbrack$ be the $k$-subalgebra of $\ol{A}$ generated by $\ol{e}$. Since $\ol{A}$ is a finite-dimensional $k$-algebra, also $\ol{B}$ is finite-dimensional. Moreover, $\ol{B}$ is commutative. Let $e \in A$ be an arbitrary element with $\phi_A(e) = \ol{e}$. Let $B \dopgleich R \lbrack e \rbrack$, a commutative subalgebra of $A$. Note that $\ol{B} = B/\fm B$. Since $A$ is a finitely generated $R$-module, the Cayley–Hamilton theorem implies that $B$ is a finitely generated $R$-algebra. Now, by the commutative case, the map $\phi_B:B \twoheadrightarrow \ol{B}$ is idempotent surjective and so there is an idempotent $e' \in B \subs A$ with $\phi_A(e') = \phi_B(e') = \ol{e}$. This shows that $\phi_A$ is idempotent surjective.
\end{proof}

The next theorem was again proven by M. Neunhöffer \cite[Proposition 6.2]{Neunhoeffer-PhD}. It is one of our key ingredients in proving Brauer reciprocity for decomposition maps in a general setting.

\begin{thm}[M.\ Neunhöffer] \label{neunhoeffer_theorem}
Suppose that $R$ is a valuation ring with fraction field $K$ and that $A$ is a finite flat $R$-algebra with split generic fiber $A^K$. If $\wh{R} \otimes_R A$ is semiperfect, where $\wh{R}$ is the completion of $R$ with respect to the topology defined by a valuation on $K$ defining $R$, then also $A$ is semiperfect.
\end{thm}

\begin{cor}[J. Müller, M. Neunhöffer] \label{algebra_dvr_semiperfect}
Suppose that $R$ is a \textit{discrete} valuation ring and that $A$ is a finite flat $R$-algebra with split generic fiber. Then $A$ is semiperfect. In particular, $\phi_A:A \twoheadrightarrow \ol{A}$ is primitive idempotent bijective and block bijective.
\end{cor}

\begin{proof}
Since $R$ is a discrete valuation ring, its valuation topology coincides with its $\fm$-adic topology so that the topological completion $\wh{R}$ is $\wh{\fm}$-adically complete, where $\fm$ denotes the maximal ideal of $R$ and $\wh{\fm}$ denotes the maximal ideal of $\wh{R}$. Hence, $\wh{R} \otimes_R A$ is semiperfect by Lemma \ref{semiperfect_std_settings}\ref{semiperfect_std_settings:complete} and Theorem \ref{neunhoeffer_semiperfect}. Now, Theorem \ref{neunhoeffer_theorem} shows that $A$ is semiperfect, too.
\end{proof}

\begin{rem}
One part of Corollary \ref{algebra_dvr_semiperfect}, the fact that idempotents lift, was also stated earlier by C.\ Curtis and I.\ Reiner \cite[Exercise 6.16]{CR-Methods-1} in an exercise in the special case where $A^K$ is assumed to be semisimple. The semisimplicity assumption was later removed by J.\ Müller in his PhD thesis \cite[Satz 3.4.1]{Muller:PhD} using the Wedderburn–Malcev theorem (this can be applied without perfectness assumption on the base field if $A^K$ splits since then $A^K/\Rad(A^K)$ is separable, see \cite[Theorem 72.19]{Curtis-Reiner-Associative-Algebras}). 
\end{rem}

\section{Further elementary facts}

Here, we prove three further elementary facts that we used in the paper. 

\begin{lem} \label{flat_faithfully_flat}
A finitely generated module $M$ over an integral domain $R$ is flat if and only if it is faithfully flat. In particular, if $M \neq 0$, we have $0 \neq \msf{k}(\fp) \otimes_R M = M(\fp)$ for all $\fp \in \Spec(R)$.
\end{lem}

\begin{proof}
We can assume that $M \neq 0$. Since $M$ is flat, it is torsion-free and so the localization map $M \rarr M_\fp$ is injective, see Lemma \ref{phi_injective_lemma}\ref{phi_injective_lemma:localiz}. Hence, $M_\fp \neq 0$. Since $M$ is a finitely generated $R$-module, also $M_\fp$ is a finitely generated $R_\fp$-module and now Nakayama's lemma implies that $0 \neq M_\fp/\fp_\fp M_\fp = \msf{k}(\fp) \otimes_R M$. Hence, $M$ is faithfully flat by \cite[Theorem 7.2]{Mat-Commutative}.
\end{proof}

\begin{lem} \label{ff_center_ext}
Let $A$ be a finite flat algebra over an integral domain $R$. Then the structure map $R \rarr A$, $r \mapsto r \cdot 1_A$, is injective. Hence, we can identify $R \subs \msf{Z}(A)$. If $R$ is noetherian, the induced map $\Upsilon: \Spec(\msf{Z}(A)) \rarr \Spec(R)$ is finite, closed, and surjective.
\end{lem}

\begin{proof}
It follows from Lemma \ref{flat_faithfully_flat} that $A$ is already faithfully flat. Let $\phi:R \rarr A$ be the structure map. This is an $R$-module map and applying $- \otimes_R A$ yields a map 
\[
A \simeq R \otimes_R A \overset{\phi \otimes_R A}{\longrightarrow} A \otimes_R A
\]
of right $A$-modules, mapping $a$ to $1 \otimes a$. This map has an obvious section mapping $a \otimes a'$ to $aa'$, hence it is injective. Since $A$ is faithfully flat, the original map $\phi$ has to be injective, too. As the image of $\phi$ is contained in the center $Z$ of $A$, the structure map is actually an injective map $R \hookrightarrow Z$. Now, assume that $R$ is noetherian. Since $A$ is a finitely generated $R$-module, also $Z$ is a finitely generated $R$-module. Hence, $R \subs Z$ is a finite ring extension and now it is an elementary fact that $\Upsilon$ is closed and surjective.
\end{proof}

The following lemma about base change of homomorphism spaces is well known but we could not find a reference in this generality (see \cite[II, \S5.3]{Bou-Commutative-Algebra-1-7} for a proof in case of a \textit{commutative} base ring).

\begin{lem} \label{base_ring_change_of_hom}
Let $A$ be an algebra over a commutative ring $R$ and let $\phi:R \rarr S$ be a morphism into a commutative ring $S$. Let $V$ and $W$ be $A$-modules. If $V$ is finitely generated and projective as an $A$-module, then there is a canonical $S$-module isomorphism.
\begin{equation}
S \otimes_R \Hom_A(V,W) \simeq \Hom_{A^S}(V^S,W^S) \;.
\end{equation}
\end{lem}

\begin{proof}
We can define a map $\gamma:S \otimes_R \Hom_A(V,W) \rarr \Hom_{A^S}(V^S,W^S)$ by mapping $s \otimes f$ with $s \in S$ and $f \in \Hom_A(V,W)$ to $s_r \otimes f$, where $s_r$ denotes right multiplication by $s$. It is a standard fact that this is an $S$-module morphism, see \cite[(2.36)]{Rei-Maximal-orders-0}). Recall that $\Hom_A(-,W)$ commutes with \textit{finite} direct sums by \cite[II, \S1.6, Corollary 1 to Proposition 6]{Bou-Algebra-1-3}. This shows that the canonical isomorphism $\Hom_A(A,W) \simeq W$ induces a canonical isomorphism $\Hom_A(A^n,W) \simeq W^n$ for any $n \in \bbN$ and now we conclude that there is a canonical isomorphism
\[
S \otimes_R \Hom_A(A^n,W) \simeq S \otimes_R W^n \simeq (S \otimes_R W)^n \simeq \Hom_{A^S}( (A^S)^n, W^S) \;,
\]
which is easily seen to be equal to $\gamma$. The assertion thus holds for finitely generated \textit{free} $A$-modules. Now, the assumption on $V$ allows us to write without loss of generality $A^n = V \oplus X$ for some $A$-module $X$. It is not hard to see that we get a commutative diagram
\[
\begin{tikzcd}[column sep=0.5ex]
S \otimes_R \Hom_A(A^n,W) \arrow{rrrrrrrrrrrrrrr}{\simeq} \arrow{d}[swap]{\simeq} &&&&&&&&&&&&&&& \left( S \otimes_R \Hom_A(V,W) \right) \arrow{d} & \oplus & \left( S \otimes_R \Hom_A(X,W) \right) \arrow{d} \\
\Hom_{A^S}( (A^S)^n, W^S) \arrow{rrrrrrrrrrrrrrr}{\simeq} &&&&&&&&&&&&&&& \left( \Hom_{A^S}(V^S,W^S) \right) & \oplus & \left( \Hom_{A^S}(X^S, W^S) \right)
\end{tikzcd}
\]
where the horizontal morphisms are obtained by the projections and the vertical morphisms are the morphisms $\gamma$ in the respective situation. The commutativity of this diagram implies that the morphism $S \otimes_R \Hom_A(V,W) \rarr \Hom_{A^S}(V^S,W^S)$ also has to be an isomorphism.
\end{proof}

\printbibliography

\end{document}